\documentclass[11pt,a4paper,reqno]{amsart}


\newtheorem{proposition}{Proposition}

\newtheorem{theorem}{Theorem}
\newtheorem{corollary}{Corollary}
\newtheorem{lemma}{Lemma}

\newtheorem{remark}{Remark}
\newenvironment{proofof}[1]{\medskip 
\noindent{\bf Proof of #1.}}{ \hfill\qed\\ }

\newcommand{\slims}{\sum\limits}

\newcommand{\To}{\longrightarrow}

\newcommand{\leqslant}{\le}

\def\eps{\varepsilon}
\def\R{{\mathbb R}}
\def\C{{\mathbb C}}
\def\N{{\mathbb N}}

\def\E{{\mathbb E}}
\def\O{{\mathcal O}}
\def\P{{\mathcal P}}
\def\1{ {\hbox{{\it 1}} \!\! I} }
\def\1{{{\mathit 1} \!\!\>\!\! I} }
\renewcommand{\phi}{\varphi}

\def\diam{{\hbox{{\rm diam}}}}
\def\Jac{{\hbox{Jac}\ }}

\def\Crit{{\mathcal C}}

\def\leb{\mbox{Leb}}
\DeclareMathOperator{\Ad}{Ad} 
\DeclareMathOperator{\Aut}{Aut} 

\begin{document}
\bibliographystyle{plain}
\title{Liv\u{s}ic regularity for Markov Systems}
\author{Henk Bruin, Mark Holland, Matthew Nicol}
\date{\today}

\begin{abstract}
We prove measurable Liv\v{s}ic theorems for dynamical systems
modelled by Markov towers. Our regularity results apply
to solutions of cohomological equations posed on  H\'enon-like
mappings  and a wide variety of  nonuniformly hyperbolic systems.
We consider both H\"{o}lder cocycles and cocycles with singularities of
prescribed order.
\end{abstract}
\maketitle
\section{Introduction}\label{sec:intro}

In this paper  we study the regularity of solutions $\psi$ of the cohomological
equation
\begin{equation}\label{eq-cohomology0}
\phi = (\psi \circ T) \psi^{-1} \qquad \mu\mbox{-a.e.}
\end{equation}
where $(T,X,\mu)$ is a dynamical system
and $\phi:X\to G$ is a cocycle taking values in a  Lie group $G$. 
Measurable rigidity in this context means that a measurable solution $\psi$
must have  a higher degree of regularity, in many contexts inheriting 
the regularity of $\phi$ and/or $T$.
Such  cohomological equations come up in different applications:
they are used to determine whether certain observables
have positive variance in the context of the Central Limit Theorem
and related distributional theorems.
In the context of group extensions, they decide on (stable) ergodicity 
and weak-mixing of the system. In other contexts, cohomological 
equations play a role in the question of whether two dynamical systems
are (H\"older or smoothly) conjugate to each other.

Fundamental work on the regularity of 
measurable solutions to cohomological  equations was
done by  Liv\v{s}ic~\cite{Livsic1,Livsic2} who established
rigidity theorems for Abelian group extensions of Anosov systems with an 
absolutely continuous invariant measure. Dynamical rigidity
theorems are often called Liv\v{s}ic theorems in the literature
because of  this.
Parry \& Pollicott~\cite{PP1}, using  a transfer operator approach,
extended Liv\v{s}ic's results to prove H\"{o}lder regularity of
coboundary and transfer functions  for compact Lie group extensions of 
subshifts of finite type and hence, via Markov partitions, Axiom
A systems. Further generalizations for uniformly hyperbolic smooth
systems are given
in~\cite{NP1,NP2,Pollicott-Walkden,Walkden1,Walkden2}. In the
Anosov setting de la Llave et al~\cite{dLMM} 
prove a $C^{\infty}$ version of   Liv\v{s}ic's  theorem and also
$C^{\infty}$  dependence of solutions upon parameters.

There are only a few results
on Liv\v{s}ic theorems for nonuniformly hyperbolic or discontinuous
systems.
 Pollicott \& Yuri~\cite{Pollicott-Yuri}  
have established Liv\v{s}ic 
theorems for H\"{o}lder $\R$-extensions of $\beta$-transformations
($f:[0,1)\To[0,1)$, $f(x)=\beta x \pmod 1$ where $\beta>1$)
via transfer operator techniques
but the regularity they obtain is  bounded variation rather than H\"{o}lder.
Jenkinson~\cite{Jenkinson} has proved that essentially bounded
measurable coboundaries  $\psi$ (i.e. solutions to $\phi=\psi\circ f-\psi$) 
for $\R$-valued smooth cocycles $\phi$ over smooth 
expanding maps $f$  have smooth versions.  

 Nicol \& Scott~\cite{NS} 
 have obtained  Liv\v{s}ic theorems for certain discontinuous  hyperbolic systems,
 showing that coboundary solutions
 taking values in Lie groups satisfying a pinching condition 
(to ensure that the system is partially hyperbolic) 
are Lipschitz if the cocycle is Lipschitz.
 The same techniques show that, for such systems,  measurable transfer functions 
taking values in compact matrix groups have Lipschitz versions. 
These results were applied  to  prove stable ergodicity for
semisimple and Abelian compact group extensions of certain uniformly hyperbolic systems
with singularities, including the $\beta$-transformation, Markov maps
and mixing Lasota-Yorke maps.

Aaronson \& Denker~\cite[Corollary 2.3]{AD}
have shown that if $(f,X,\mu,\P)$ is a mixing Gibbs-Markov
map preserving a probability measure $\mu$ with countable
Markov partition $\P$ and $\phi: X\to \R^d$ is Lipschitz
(i.e. $d(\phi(x), \phi(y)) \leq C \rho(x,y)$ for all $x,y \in f(a)$,
$a\in \P$, with respect to a metrics $d$ on $\R^d$ and $\rho$ on $X$ 
derived from the symbolic dynamics), 
then any measurable solution $\phi:X\to \R^d$
to $\phi=\psi\circ f -\psi$ has a version $\tilde \phi$
which is Lipschitz continuous.

The work of \cite{ADSZ} is a study of the statistical properties of fibred
systems and gives   rigidity results which provide checkable conditions for
the aperiodicity of cocycles (i.e. nonexistence
of solutions $\psi$) which allow one to establish, for example,
 de Moivre's approximation for various systems, including the
 $\beta$-transformation. A related  result is given in ~\cite[Lemma 6.1.2]{Gou}.

In two influential papers ~\cite{Young1,Young2}, Young 
describes properties of a class of Markov extensions (which
we will call  Young towers) which are useful  
to establish  rates of decay of correlations
and the CLT in non-uniformly hyperbolic systems. Scott~\cite{Scott,thesis} has recently proved 
measurable Liv\v{s}ic theorems for 
certain Lie group valued H\"{o}lder cocycles over a class of unimodal maps
modelled by a Young tower~\cite{Young1}. 
More precisely suppose
$(f,X,\mu)$  ($X$ a finite collection of  intervals)
 is a unimodal map (belonging to a certain class) and
$g:X\to G$ is a  Lie group valued H\"{o}lder cocycle
(satisfying a pinching condition if $G$ is noncompact). If
$\psi$ is a  measurable solution to $g=(\psi\circ f) \psi^{-1}$ $\mu$-a.e.
then $\psi$ is H\"{o}lder on an arbitrarily large open set (i.e.
given $\eps>0$ there exists an open set $U$ such that 
$\psi$ is H\"{o}lder on $U$ and $\mu(U)>1-\eps$). Similar
measurable Liv\v{s}ic theorems for other types of cohomological equations are given.

This paper extends the results of Scott in several directions. We prove
measurable Liv\v{s}ic theorems for 
more general Markov extensions and for cocycles with singularities of prescribed
order. We also obtain  regularity results for
measurable  conjugacies between certain non-uniformly hyperbolic
systems and measurable Liv\v{s}ic theorems for certain  non-uniformly
hyperbolic systems (in particular H\'enon-like mappings).

After writing this paper  we learnt that Gou\"ezel \cite{Gou3} has obtained 
similar results for cocycles into Abelian groups over one-dimensional 
Gibbs-Markov systems and Young towers. 
From  ~\cite{Gou3} we learnt
of  a Martingale Density  Theorem (see Appendix) which allows a  more elegant
approach in part of our proof than  the argument 
using Lebesgue density points adapted from~\cite{NS}.

{\em Acknowledgements:} We would like to thank Ian Melbourne and Gerhard 
Keller for fruitful suggestions. M. Holland also acknowledges the 
support of an EPSRC grant, No. GR/S11862/01. M. Nicol acknowledges the 
support of NSF grant DMS-0071735.

\section{Cohomological equations and group extensions}\label{sec:groups}

In this section we collect some facts about group extensions and
cohomological equations that 
set the framework of  this paper. All Lie groups are assumed connected
and  finite-dimensional.

Suppose $f: X\to X$ is an ergodic dynamical system with
respect to an invariant measure $\mu$. Let $G$ be a compact  Lie
group and let $dh$ denote Haar measure on $G$. Suppose $\phi:
X\to G$ is measurable.
\begin{proposition}\cite{Keynes-Newton,noorani}
\begin{enumerate}
\item
The compact group extension $T(x,g)=(f(x),\phi (x)g)$ is ergodic
with respect to $\mu \times dh$ if and only if the equation 
\begin{equation}\label{eq-cohomology1}
\psi(fx)=R(\phi (x))\psi(x) \qquad \mu\mbox{-a.e.}
\end{equation}
where $R$ is an irreducible (unitary) representation of degree $d$
and $\psi: X \To \C^d$ is measurable, is only satisfied when 
$\psi$ is constant or $R$ is the trivial representation. 
\item Suppose  $f: X\to X$
is weak-mixing and $T: X\times G\to X\times G$ is ergodic.
Then $T(x,g)=(f(x),\phi (x)g)$ is weak-mixing  with respect to $\mu
\times dh$ if and only if for any $e^{i\alpha}\not = 1$ and any
non-trivial one-dimensional representation $\chi$ of $G$ the
equation
\begin{equation}\label{eq-cohomology2}
 \psi(fx)=e^{i\alpha} \chi(\phi (x))\psi(x) \mbox{  $\mu$-a.e. }
\end{equation}
has no nontrivial measurable solution $\psi: X\to \C$.
\end{enumerate}
\end{proposition}

Note that the aperiodicity condition of~\cite{ADSZ} 
\[
\gamma\circ \phi=\lambda \psi / \psi\circ f
\]
 where $\gamma$ is a character
of $G$ (a locally compact Abelian polish group), $\lambda \in S^1$
is a special case of equation~\eqref{eq-cohomology2}.

Suppose that $\phi_i:X\to G$, $i = 1,2$ are two compact Lie-group valued 
cocycles over a system $(f,X,\mu)$. A measurable function
$\psi:X \to G$ conjugates the $G$ extensions $T(x,g)=(fx,\phi_i(x) g)$ ($i=1,2)$ if 
\begin{equation}\label{eq-cohomology3}
\psi(fx)\phi_1(x)=\phi_2(x)\psi (x)\qquad \mu\mbox{-a.e.}
\end{equation}
We call such a conjugating function $\psi$ a {\em transfer function}.

If $G$ is compact we may identify $G$ with a subgroup of $U(d)$, the group of 
$d\times d$ unitary matrices.
In  this representation we may identify a $G$-valued $\phi:X\to G$ with 
$\phi: X\to \C^{d^2}$. Define $\theta(x):M(d)\to M(d)$, 
 a mapping from the space of $d\times d$ complex matrices to itself  by
\[
\theta(x): A \to \phi_2(x) A \phi_1(x)^*.
\]
It is possible to show that $\theta(x): M(d)\to M(d)$ is unitary. 
There is a standard way, see~\cite[Theorem 1]{PP1} and \cite[Theorem A]{Jenkinson} to rewrite 
$\psi(fx)\phi_1(x)=\phi_2(x)\psi (x)$
in form $\psi(fx)=\theta (x) \psi (x)$ where $\psi:X\to \C^{d^2}$, $\theta(x)\in U(d)$.
Hence  the question of the regularity of conjugacies between compact group extensions 
may be reduced to those of 
the regularity of solutions to  equation (\ref{eq-cohomology1}).

The proof of  our  coboundary  Liv\v{s}ic regularity results, such as 
Theorem~\ref{main}, may be slightly modified (as in \cite[Section 2.1]{NS}), to establish the 
same degree of regularity for solutions $\psi$ to 
 equation~\eqref{eq-cohomology1}, equation~\eqref{eq-cohomology2} 
 or equation~\eqref{eq-cohomology3} posed over the same dynamical
 system. We omit the straightforward  proof of this and refer the reader to 
 \cite[Section 2.1]{NS}.

\subsection{Lie groups}

Let $G$ be a connected Lie group with Lie algebra denoted by
$L(G)$ which we identify with the tangent space at the identity,
$T_eG.$ We let $r_g$ denote right multiplication by $g\in G$.
Given a norm $\|\cdot\|$ on $T_eG$ we define a norm on $T_gG$ by
$\|\mathbf{v}\|_g=\|r_{g^{-1}}\mathbf{v}\|_e$. This norm induces a
right invariant metric $d_G$ on  $G$ so that
$d_G(gk,hk)=d_G(g,h)$, see \cite[Section 4]{Pollicott-Walkden}.
Throughout this paper we will write $d(\cdot,\cdot)$ instead of
$d_G(\cdot,\cdot)$ when it is clear from context that we mean the
metric on $G$. For a general reference on Lie groups 
see~\cite{BrockDieck}.

We define the adjoint map $\Ad\colon G\to\Aut(L(G))$, for $g\in G$
and $X\in T_eG$ by
\[
\Ad(g)\mathbf{v}=\frac{d}{dt}(g \exp(t\mathbf{v})g^{-1})
\mbox{ at } t=0.
\]
Note that when $G$ is a matrix group this action is conjugation i.e. 
$v\to gvg^{-1}$.
A calculation ~\cite[Section 4]{Pollicott-Walkden} shows that
\begin{equation}\label{eq-Gmetric_left_property}
d(gh,gk)\leq \|\Ad(g)\|d(h,k).
\end{equation}

Suppose $\phi\colon M\to G$ is H\"{o}lder of exponent $\alpha>0.$ 
Define
\begin{eqnarray*}
\mu_u&:=&\lim\limits_{n\rightarrow\infty}\left(\sup\limits_{x\in
X}\|\Ad(\phi_n(x))\|\right)^{\frac{1}{n}},
\end{eqnarray*}
where $\phi_n (x) =\phi(f^{n-1}x)\ldots \phi(fx)\phi(x)$. 

If $G$ is Abelian or compact, then $\sup_x \|\Ad(\phi_n(x))\|$ is bounded in $n$, 
whereas if $G$ is nilpotent,
$\sup_x\|\Ad(\phi_n(x))\|$ can grow at most at a polynomial rate.
In these three cases, $\mu_u=1$. 

\section{Axiomatic Approach for Nonuniformly Expanding Maps}\label{sec:markov}

Let $(M,\rho,\mu)$ be a metric space endowed with a non-atomic Borel 
probability measure $\mu$. We assume $M$ can be decomposed as
$M = \cup_k M_k \bmod \mu$, where each $M_k$ is connected and
$\sup_k (\diam (M_k)) \leq 1$.
Let $f:\cup_k M_k \to M$ be a map  such that $f|M_k$ is continuous for
each $k$ and 
such that $\mu$ is $f$-invariant and ergodic.

Let $\P_0$ be the partition of $M$ into the sets $M_k$, and
$\P_n = \bigvee_{i = 0}^{n-1} f^{-i}(\P_0)$.
For $x \in M$, let $\P_n[x]$ be the partition element (cylinder set) 
in $\P_n$ containing $x$.

Consider the natural extension $(\hat M, \hat f, \hat \mu)$
of $(M,f,\mu)$: each point $\hat x \in \hat M$ is a sequence
\[
\hat x = (x_0, x_1, x_2, \dots ) \mbox{ with } M \owns x_i = f(x_{i+1})
\mbox{ for all } i \geq 0.
\]
The measure $\hat \mu$ is defined in the standard way ~\cite{Keller}
and in particular for each $n$:
\[
\hat \mu(\{ \hat x \in \hat M \ | \ x_n \in A \}) = \mu(A)
\]
for each $\mu$-measurable set $A$.

We assume:
\begin{enumerate}
\item For all $k$, 
$f:M_k \to f(M_k)$ is one-to-one and $f(M_k)$ is equal to a  union 
of components $M_l$ $\bmod~\mu$ (Markov property).
\item There exists $\lambda>1$ and for $\hat \mu$-a.e. $\hat x$ a number  
$K(\hat x)$ such that
\begin{equation}\label{expo_contraction}
\rho(y_n,z_n)\leq K (\hat x)\lambda^{-n} \rho(y_0,z_0)
\end{equation}
for all $n \geq 0$ and $y_n,z_n  \in \P_n[x_n]$.
(In dimension one this assumption can be weakened, see Section~(\ref{sec_onedim})).
\item Let $J_{\mu}(x)$ denote the Jacobian of $\mu$ at $x$.
For $\hat x \in \hat M$, define $J_{\mu}^n(x_n) = 
\prod_{i=1}^n J_{\mu}(x_i)$.
For $\hat \mu$-a.e. $\hat x$, there exists a constant $C(\hat x)$ such that
if $\hat y, \hat z \in  \hat M$ are such that $y_i, z_i \in \P_i[x_i]$ 
for all $i$, then
\begin{eqnarray}\label{BD}
\left| \frac{ J_{\mu}^n(y_n) }{ J_{\mu}^n(z_n) } \right|\le  C(\hat x) .
\end{eqnarray}

\end{enumerate}

The Jacobian $J_{\mu}$ of $\mu$ is $\gamma$-H\"older
if there exists $C$ and $\gamma \in (0,1]$ such that
\begin{equation}\label{Jac_Holder}
\left| \frac{J_\mu(x)}{J_{\mu}(y)} - 1 \right| \leq 
C \cdot \rho(f(x),f(y))^{\gamma}.
\end{equation}
H\"olderness of the Jacobian implies a result stronger than  (\ref{BD}).
\begin{lemma}\label{dist_Jac}
Assume $J_{\mu}$ is $\gamma$-H\"{o}lder with coefficient $C$. 
For $\hat \mu$-a.e. $\hat x$, there exists a constant 
$B = B(\hat x,C,\lambda^{\gamma})$ such that
\[
\left| \frac{ J_{\mu}^n(y_n) }{ J_{\mu}^n(z_n) } \right|\le
1 + B \rho(y_0 , z_0)^{\gamma}.
\]
\end{lemma}

\begin{proof}
Using \eqref{Jac_Holder} and \eqref{expo_contraction}, we obtain 
\begin{eqnarray*}
\left| \frac{ J_{\mu}^n(y) }{ J_{\mu}^n(z) } \right| &=&
\prod_{i=0}^{n-1} \frac{|J_{\mu}(y_i)|}{|J_{\mu}(z_i)|}  \\
&\leq& \prod_{i=0}^{n-1} ( 1 + C \rho(f(y_i) , f(z_i))^{\gamma}) \\
&\leq& \exp( C \sum_{i=1}^{n} \rho(y_i,z_i)^{\gamma} ) \\
&\leq& \exp( C \cdot K(\hat x) \sum_{i=1}^{n} \lambda^{-i \gamma} 
\rho(y_0,z_0)^{\gamma} ) \\
&\leq& \exp\left( \frac{C \cdot K(\hat x)}{\lambda^{\gamma} - 1} 
\rho(y_0 ,z_0)^{\gamma} \right), \\
\end{eqnarray*} 
which is smaller than $1+ B \rho(y_0 ,z_0)^{\gamma}$
for some $B$ depending only on $\hat x$, $C$, $\lambda^{\gamma}$
and the diameter of the component $M_k$ containing $x_0$.

In fact, the above computation only requires that
$\sum_{i=0}^{n-1} \diam(\P_i[x_i])^\gamma \leq K_0 (\hat x) < \infty$
for $\mu$-a.e  $\hat x $ independently of $n$, which is an estimate 
valid under a less strict assumption than (\ref{expo_contraction}).
\end{proof}

\noindent {\bf Cocycle assumptions:}
Let $\phi\colon M\to G$ be H\"{o}lder of exponent $\alpha>0$. 
Recall that
\begin{eqnarray*}
\mu_u&:=&\lim\limits_{n\rightarrow\infty}\left(\sup\limits_{x\in
M}\|\Ad(\phi_n(x))\|\right)^{\frac{1}{n}},
\end{eqnarray*}
where $\phi_n(x)=\phi(f^{n-1}x)\ldots \phi(fx)\phi(x)$ with $\phi_0(x) = e$.
If $G$ is Abelian, $\| \Ad(\phi_n(x)) \| = 1$ 
and $\| \Ad(\phi_n(x)) \|$ is bounded
if $G$ is compact. For  nilpotent groups $G$,
$\| \Ad(\phi_n(x)) \|$ grows at most at a polynomial rate in $n$, 
so $\mu_u = 1$.
For the general case, we impose a partial hyperbolicity condition (\ref{PH})
on the group extension:
\begin{equation}\label{PH}
 1\leq  \mu_u < \lambda^{\alpha} \tag{PH}
\end{equation}
where $\lambda$ is from (\ref{expo_contraction}).

\begin{theorem}\label{main}
Assume that $(M,f,\mu)$ is a measure preserving  
Markov system as above and let $M_k \in \P_0$.
Let $\phi:M \to G$ be a Lie group valued $\alpha$-H\"older observable 
(i.e. $d(\phi(x),\phi(y)) \leq C \rho(x,y)^{\alpha}$) satisfying
the partial hyperbolicity condition (\ref{PH}) above. 
Let  $\psi:M \to G$ be a $\mu$-measurable solution of the cohomological 
equation
\[
\psi \circ f(x) = \phi(x) \cdot \psi(x) \qquad \mu\mbox{-a.e.}
\]
Then there is a version $\tilde \psi$ of $\psi$ (i.e. $\psi = \tilde \psi$ 
$\mu$-a.e.) such that $\tilde \psi$ is $\alpha$-H\"older on $M_k$. 
\end{theorem}

\begin{corollary}\label{cor_extension}
If $f^j (\cup_{k \in S} M_k)=M$ for some $j>0$ and finite collection
of indices $S$ then there is a version which is $\alpha$-H\"{o}lder
on $M$.
\end{corollary}

\begin{proofof}{Corollary~\ref{cor_extension}} 
By considering the cohomological equation 
$\psi\circ f^n (x)=\phi(f^{n-1} x)\cdot
\phi(f^{n-1} x) \ldots \phi (x) \cdot \psi (x)$
we may extend the version of $\psi$ as a H\"{o}lder
function to any image $f^j(M_k)$.
\end{proofof}

\begin{remark}\label{RemForward}
It is easy to show that given
$\eps>0$ there is a version of $\psi$ which is $\alpha$-H\"{o}lder
on a finite union of sets $\cup_{k \in S} M_k$ such that
$\mu(\cup_{k \in S} M_k)>1-\eps$.
The H\"{o}lder coefficient depends in general upon $S$  but the exponent is uniform.
\end{remark}

\begin{remark}
A slight modification of  the proof
shows that the same regularity results  hold for  solutions $\psi$ to 
equation~\eqref{eq-cohomology1},
equation~\eqref{eq-cohomology2} or  equation~\eqref{eq-cohomology3}.
\end{remark}

\begin{proofof}{Theorem~\ref{main}}
Choose any $\Lambda := M_k \in \P_0$ such that $\mu(\Lambda) > 0$.
Let $0<\delta<1$. As a consequence of the Martingale Density Theorem 
(see Appendix)  
for $\hat \mu$-a.e. $\hat x \in \hat M$  and for 
infinitely many $n:$
\[
 \frac{ \mu \{ y_n \in \P_n [x_n]: 
d( \psi(y_n), \psi(x_n))<\delta \} }{ \mu(\P_n [x_n])} 
> 1-\delta.
\]
Let $\hat x$ be such a point with $x_0\in \Lambda$.
For simplicity of notation in the rest of the proof we will not indicate the dependence of constants
upon $\hat x$.
We consider points
$\hat y=(y_0,y_1,\ldots y_n\ldots )$ and $\hat z=(z_0,z_1,\dots z_n\dots)\in\hat M$ 
such that $y_n,z_n \in \P_k[x_n]$ for all $n=0,1,\ldots$. 
Hence $\hat x$, $\hat y$ and $\hat z$ are all paths in the ``same inverse branch'' of 
$f$.

On $\Lambda$ we define a function $\Phi:\Lambda \to G$ by
\[
\Phi(y_0)=\lim_{n\to\infty}\phi_n(y_n)\phi_n(x_n)^{-1},
\]
where $\phi_n(x_n)=\phi(x_1)\ldots \phi(x_n).$ 
This function is well defined since
\begin{eqnarray*}
\lefteqn{d(\phi_{n+1}(y_{n+1})\phi_{n+1}(x_{n+1})^{-1},\phi_n(y_n)\phi_n(x_n)^{-1})}\\
&=&d(\phi_n(y_n)\phi(y_{n+1})\phi(x_{n+1})^{-1}\phi_n(x_n)^{-1},\phi_n(y_n)\phi_n(x_n)^{-1})\\
&=&d(\phi_n(y_n)\phi(y_{n+1})\phi(x_{n+1})^{-1}\phi_n(x_n)^{-1},\phi_n(y_n)\phi(x_{n+1})
\phi(x_{n+1})^{-1}\phi_n(x_n)^{-1})\\
&=&d(\phi_n(y_n)\phi(y_{n+1}),\phi_n(y_n)\phi(x_{n+1}))\quad\quad\mbox{
(right invariance)}\\
&\leq&\|\Ad(\phi_n(y_n))\| \ d(\phi(y_{n+1}),\phi (x_{n+1}))\quad\quad\mbox{(by
\eqref{eq-Gmetric_left_property})}\\
&\leq&C \cdot K(\hat x) \left((\mu_u )\lambda^{-\alpha}\right)^n
\quad\quad\mbox{(by
\eqref{expo_contraction})}\\
&\leq&C \cdot K(\hat x) \cdot \kappa^n,
\end{eqnarray*}
where $\kappa \in (0,1)$ by (PH).
Thus the sequence $\phi_n(y_n)\phi_n(x_n)^{-1}$ is Cauchy and so
converges.

Next we show that $\Phi$ is H\"{o}lder. Let $y_0,z_0\in \Lambda$, then
\begin{eqnarray}\label{eq_Holder}
\lefteqn{ d(\phi_n(y_n)\phi_n(x_n)^{-1},\phi_n(z_n)\phi_n(x_n)^{-1}) } 
\notag \\
&=& d(\phi_n(y_n),\phi_n(z_n))  \\
&\leq&\slims_{i=0}^{n-1}d(\phi_i(y_i)\phi(y_{i+1})\phi_{n-i-1}(z_{i+1}),
\phi_i(y_i)\phi(z_{i+1})\phi_{n-i-1}(z_{i+1})) \notag\\
&\leq&\slims_{i=0}^{n-1}\|\Ad(\phi_i(y_i))\| \ d(\phi(y_{i+1}),\phi(z_{i+1})) \notag\\
&\leq&\slims_{i=0}^{n-1}C \cdot K(\hat x) (\mu_u  )^{i+1} \lambda^{-(i+1)\alpha}\rho(y_0,z_0)^{\alpha}. \notag
\end{eqnarray}
Letting $n\rightarrow\infty$ gives
$d(\Phi (y_0),\Phi (z_0))\leq C \cdot K(\hat x) \cdot \rho(y_0,z_0)^{\alpha}$.
It is clear that if  $\phi$ is  Lipschitz (i.e. $\alpha=1$)
then  $\Phi$ is also Lipschitz.

Define
\[
\Psi_n(y_0)=\phi_n(y_n)\phi_n(x_n)^{-1}.
\]
Then 
\begin{eqnarray*}
\psi(y_0)&=&\phi_n(y_n)\psi(y_n)\\
&=&\Psi_n(y_0)\phi_n(x_n)\psi(x_n)\psi(x_n)^{-1}\psi(y_n)\\
&=&\Psi_n(y_0)\psi(x_0)\psi(x_n)^{-1}\psi(y_n).
\end{eqnarray*}
Thus
\begin{eqnarray*}
d(\psi(y_0),\psi(z_0))&\le& d(\Psi_n(y_0)\psi(x_0)\psi(x_n)^{-1}\psi(y_n),
\Psi_n(y_0)\psi(x_0))\\
&+&d(\Psi_n(y_0)\psi(x_0)\psi(x_n)^{-1}\psi(z_n),\Psi_n(y_0)\psi(x_0))\\
&+& d(\Psi_n(y_0)\psi(x_0)\psi(x_n)^{-1}\psi(z_n),
\Psi_n(z_0)\psi(x_0)\psi(x_n)^{-1}\psi(z_n)).
\end{eqnarray*}

By right-invariance of the metric the last term may be written as
\[
d(\Psi_n(y_0), \Psi_n(z_0)).
\]
As a function  of $y_0$, $\Psi_n(y_0)$ converges to the
$\alpha$-H\"{o}lder function $\Psi(y_0)$. Thus letting 
$n\rightarrow \infty$, we obtain $d(\Psi(y_0), \Psi(z_0))\le C \rho(y_0,z_0)^{\alpha}$.

Given $\eta>0$ there exists $\delta_{\eta}>0$
such that $d(\psi(z_n),\psi(y_n))\le \delta_{\eta}$ implies 
\begin{eqnarray*}  
 d(\Psi_n(y_0)\psi(x_0)\psi(x_n)^{-1}\psi(y_n),
\Psi_n(y_0)\psi(x_0))&\le& \frac{\eta}{2},\\
d(\Psi_n(y_0)\psi(x_0)\psi(x_n)^{-1}\psi(z_n),\Psi_n(y_0)\psi(x_0))
&\le&\frac{\eta}{2}.
\end{eqnarray*}

Choose $n$ sufficiently large so that 
\[
 \frac{ \mu \{ y_n \in \P_n [x_n]: 
d(\psi(y_n),\psi(x_n))<\delta_{\eta} \} }{ \mu(\P_n [x_n])} 
> 1-\delta_{\eta}.
\]
Now we estimate $\mu(f^n( \{ y_n \in \P_n [x_n]: 
d(\psi(y_n),\psi(x_n))<\delta_{\eta} \}))$ relative to $\mu(\Lambda)$.
By boundedness of distortion of the Jacobian of $f^n$ we have that
\begin{eqnarray}\label{Hmc}
&& \frac{\mu(f^n(\{ y_n \in \P_n [x_n]: 
d(\psi(y_n),\psi(x_n))<\delta_{\eta} \}))}{\mu(f^n\P_n[x_n])} \notag \\
&& \qquad \qquad
\leq \O(1)\ \frac{\mu(\{ y_n \in \P_n [x_n]: 
d(\psi(y_n),\psi(x_n))<\delta_{\eta} \})}{\mu(\P_n [x_n])}.
\end{eqnarray}
Hence for the above $\eta>0$, choosing $\delta_{\eta}$ smaller if necessary, we have 
$$\mu(f^n(\{ y_n \in \P_n [x_n]: 
d(\psi(y_n),\psi(x_n))<\delta_{\eta} \})) > (1-\eta)\mu(\Lambda).$$
Since  
$d(\psi(y_n),\psi(z_n))<\delta_{\eta}$  implies 
$d(\psi(z_0),\psi(y_0))<\eta +K \rho(z_0, y_0)^{\alpha}$ we have shown that
\begin{equation*}
\begin{split}
\mu\times\mu &\{(y_0,z_0)\in \Lambda \times \Lambda:
d(\psi(z_0),\psi(y_0))<2\eta+2K \rho(z_0,y_0)^{\alpha}\}\\
 &\qquad\qquad\qquad\qquad\qquad\qquad\qquad\qquad
  >(1-2\eta)\mu\times\mu(\Lambda \times \Lambda).\\
\end{split}
\end{equation*}
Since $\eta$ was arbitrary, $\psi\mid \Lambda$ has a H\"older version.
\end{proofof}

\subsection{Cocycles with singularities}\label{sec:sing}

Let $\varphi:M\to G$ be a cocycle which is H\"older 
except for discontinuities and singularities concentrated on a finite set $\mathcal{C}$. 
Let $\lambda$ be as in (\ref{expo_contraction}), 
which is defined $\mu$-a.e, and for a fixed $\delta>0$
let $B(c,\delta)$ denote a $\delta$-neighbourhood of $c$ for $c\in\mathcal{C}$.
We consider the following three scenarios:
\begin{enumerate}
\item {\bf Bounded discontinuity:} The cocycle $\varphi(x)$ is bounded and $\gamma$-H\"older in 
the complement of $\mathcal{C}$ but for each $c\in\mathcal{C}$ we have: \\
$\lim_{x\to c^{+}}\varphi(x)\neq\lim_{x\to c^{-}}\varphi(x)$.

\item {\bf Logarithmic singularity:} For each $c\in\mathcal{C}$ inside $B(c;\delta)$ we have 
$\| \varphi(x)\|\approx |\log \rho(x,c)|$, where $\| \ \|$ denotes
the norm of the group element. Moreover there is
a sequence $\{\varepsilon_n\},$ such that 
$\sum_n \mu(B(c;\varepsilon_n )) < \infty,$ and 
 $$\limsup_{n\to\infty}\frac{\log\varepsilon^{-1}_{n}}{n\log\lambda}<1.$$ 
\item {\bf Pole:} For each $c\in\mathcal{C}$ there is  $p>0$ so that on $B(c;\delta)$, 
$\|\varphi(x)\|\approx \rho(x,c)^{-p}$, and for some  
sequence $\{\varepsilon_n\},$ such that $\sum_n \mu(B(c;\varepsilon_n )) < \infty,$ 
we have
$$\limsup_{n\to\infty}\frac{(p+1)\log\varepsilon^{-1}_{n}}{n\log\lambda}<1.$$
\end{enumerate}
In each case, we will assume that the Holder exponent of $\varphi$ restricted to the 
complement of $\cup_{c\in\mathcal{C}}B(c;\delta)$ is $\gamma$. 
We state the following result:
\begin{theorem}\label{nonholder}
Assume that $(M,f,\mu)$ is a measure preserving smooth Markov system as 
defined in Section~\ref{sec:markov},
and $M_k \in \P_0$. Let $\varphi:M\to G$ be a Lie group 
valued observable which has a singularity set 
$\mathcal{C} \subset M \setminus \cup_k \overline{M_k}$, 
characterized by cases either (1), (2) or (3) above. 
Let $\psi$ be a $\mu$-measurable solution of the cohomological equation
$$\psi\circ f(x)=\varphi(x)\psi(x)\qquad \mu\mbox{-a.e.}$$ 
Then there is a version $\tilde\psi$ of $\psi$, with 
$\tilde\psi=\psi$ $\mu$-a.e 
such that $\tilde\psi$ is $\alpha$-H\"older on $M_k$, 
for some $\alpha\in(0,1)$.  
\end{theorem}

\begin{remark}
The condition that $\mathcal{C} \subset M \setminus \cup_k \overline{M_k}$
is unnecessarily strong, but is shared by many examples.
For multimodal maps and Lorenz maps on the interval,
induced maps $F$ to a neighborhood $Y$ of the critical point $c$ are common
constructions ~\cite{BLS}. The interval $Y$ has the decomposition  
$Y= \cup_i Y_i \bmod \mu$
where $F|Y_i = f^{r_i}:Y_i \to Y$
is monotone onto for an appropriate $r_i > 0$. For each $i$,
$\overline{f^j(Y_i)} \not\owns c$ for $0 \leq j < r_i$.
Hence, if $\varphi$ has only singularities at $c$, then
using the above theorem and by the argument of  Corollary~\ref{cor_extension}, we can conclude that
$\psi$ has a H\"older version on $Y$.
\end{remark}

\begin{remark}
Also in cases where $\mathcal{C} \subset M \setminus \cup_k \overline{M_k}$
fails, the proof below can still be used to get partial results.
Given an element $M_k$ the proof constructs an $N$ and a component $J$ (where $f^N(J)=M_k$) 
of the preimage $f^{-N}(M_k)$  such that $\psi$ has a H\"older version
on $J$. If $J$ can be chosen such that 
$\overline{f^i(J)} \cap \mathcal{C} = \emptyset$
for $0 \leq i < N$, then there is a version of $\psi$ which is H\"older on $M_k$.
If this condition is not satisfied then using the proof of Corollary~\ref{cor_extension}
it is possible to show that for any $\eps>0,$  $M_k$ 
contains an open set $U \subset M_k$,
$\mu (M_k \setminus U)<\eps$,
and $\psi$ has a version which is H\"{o}lder on $U$.
\end{remark}

\begin{remark} The dependence  of $\alpha$ on the exponent $\gamma$, the 
asymptotics of the sequence $\varepsilon_n$, and the type of singularity
are apparent from  the proof. 
\end{remark}
 
\begin{remark}
The same regularity is forced upon solutions $\psi$ to equation~\eqref{eq-cohomology1},
equation~\eqref{eq-cohomology2} or  equation~\eqref{eq-cohomology3}.
\end{remark}

\begin{proof}
As in the proof of Theorem~\ref{main}, choose any $\Lambda=M_k \in
\P_0$ such that
$\mu(\Lambda)>0$ and let $(\hat M, \hat f, \hat \mu)$ denote the 
natural extension of $(M,f,\mu)$.
We  have to check that there are sufficiently many backward paths that
avoid passing too close to the singularity. We do this by using a Borel-Cantelli argument.
Let $B_n=B(c;\varepsilon_n)$, then $\sum_{n\geq 1}\mu(B_n)<\infty$, and hence we deduce that
for  $\hat \mu$-a.e. $\hat x \in \hat M$,
there exists $N = N(\hat x)$ such that $x_n \notin B_n$ for all $n \geq N$.
Combining all these facts we obtain that for $\hat\mu$-a.e.  
backward orbit there exists $N(\hat x)$ such that for all $n \geq N(\hat x)$
the following hold simultaneously:
\begin{eqnarray}
&& x_0 \in \Lambda \mbox{ and } f^n(\P_n[x_n]) = \Lambda \label{onto}. \\
&& \diam( f^{-n}(\Lambda) \cap \P_n[x_n] ) \leq \frac{1}{C(\hat x)} \lambda^{-n}. \label{expo} \\
&& \rho( x_n , c ) \geq 2 \varepsilon_n.  \label{holder}
\end{eqnarray}

The last two observations show that $\varphi\mid f^{-n}(\Lambda)\cap 
\P_n[x_{n}]$ is a H\"older function for $n\ge N$.
We now consider  the H\"older properties of $\varphi$ in the cases
that we are interested in. 
Suppose for $n \geq N$ we have
$y,z \in  f^{-n}(\Lambda) \cap \P_n[x_{n}]$.
In the case of a logarithmic singularity we have (inside $B(c;\delta)$):
\begin{eqnarray*}
d(\varphi(y),\varphi(z)) &\approx& |\log \rho(y,c) - \log \rho(z,c)|
\leq\log\biggl(1 + \frac{\rho(y,z)}{\rho(z,c)}\biggr) \\
&\leq& \rho(y,z)^{1-\alpha} \frac{(\lambda^{-n}/C(\hat x))^{\alpha}}{\rho(z,c)}
\leq \frac{C(\hat x)^{-\alpha}}{\varepsilon_{n}\lambda^{n\alpha}} 
\rho(y,z)^{1-\alpha}\\
\end{eqnarray*}
for some constant $\alpha>0$. For $\iota>0$  arbitrary, we then obtain the estimate
\[
d(\varphi(y),\varphi(z))\leq \tilde C (\hat x)
\rho(y,z)^{1-\tilde\alpha-\iota}\quad\textrm{with}
\quad\tilde\alpha=\limsup_{n\to\infty}
\frac{\log\varepsilon^{-1}_{n}}{n\log\lambda}<1.
\]
Outside $B(c;\delta)$, the function $\varphi$ will be $\gamma$-H\"older.
 
Now consider the case where $\varphi$ has a finite order pole. 
Arguing as in the case of a logarithmic singularity we obtain (inside $B(c,\delta)$):
\begin{eqnarray*}
d(\varphi(y),\varphi(z)) &\approx& |\rho(y,c)^{-p} - \rho(z,c)^{-p}|\\
&\leq& \max \bigl\{ p \rho(z,c)^{-p-1},p \rho(y,c)^{-p-1}\bigr\}\rho(y,z)\\
&\leq&\frac{C(\hat x)^{-\alpha}}{(\varepsilon_{n})^{p+1}\lambda^{n\alpha}} 
\rho(y,z)^{1-\alpha}\\
\end{eqnarray*}
and hence for arbitrary $\iota>0$ we obtain
$$d(\varphi(y),\varphi(z))\leq \tilde{C}(\hat x)
\rho(y,z)^{1-\tilde\alpha-\iota}\quad\textrm{with}
\quad\tilde\alpha=\limsup_{n\to\infty}\frac{(p+1)\log\varepsilon^{-1}_{n}}{n\log\lambda}<1.$$
Outside $B(c;\delta)$, $\varphi$ will be $\gamma$-H\"older.

In the case of a bounded discontinuity, the Borel-Cantelli argument
is simpler, since we only have to worry about 
$f^{-n}(\Lambda) \cap P_n[x_n]$ intersecting $\mathcal{C}$,
which is impossible by the assumption that $\mathcal{C} \cap \overline M_k
= \emptyset$.

So we proved now that $\psi$ has a H\"older version on
$f^{-N}(\Lambda) \cap P_N[x_N]$.
To show that $\psi\mid \Lambda$ has a H\"older version, we argue as 
Corollary~\ref{cor_extension}. The fact that $\mathcal C$ is disjoint
from each $\overline M_k$  implies
that the version $\tilde\psi$ will be H\"{o}lder.
\end{proof}

\section{One-Dimensional Systems}\label{sec_onedim}

In this section, we consider $C^{2}$ one-dimensional systems 
for which $\rho$ is Euclidean distance and
$\mu$ is an invariant measure which is absolutely continuous with 
respect to Lebesgue. The assumption that there exists a
function  $K(\hat x)$ and $\lambda>0$ such that
\[
\rho(y_n,z_n)\leq K (\hat x)\lambda^{-n}\rho(y_0,z_0)
\]
for all $n \geq 0$ and $y_n,z_n\in \P_n[x_n]$ can be replaced by
two conditions which are commonly assumed in the literature,
\begin{enumerate}
\item $f$ is $C^2$  and has bounded distortion uniformly over all iterates:
there exists a function  $K(\hat x)$ such that
\begin{equation}\label{distortion}
\frac{|Df^n(y)|}{|Df^n(z)|} \leq K (\hat x)
\end{equation}
for all $n \geq 0$ and $y,z  \in \P_n[x_n]$.
\item Positive  Lyapunov exponent, i.e.  
$\lambda(\mu) = \exp \int \log |Df| d\mu > 1$.
\end{enumerate}

Instead of  assuming (\ref{expo_contraction}) we may use (\ref{distortion}) and 
$\lambda(\mu)>1$  to prove the following lemma.

\begin{lemma}\label{lem_expo_contraction}
For $\hat \mu$-a.e. $\hat x \in \hat M$, there exists a constant
$C(\hat x)$ such that
\begin{equation}\label{expo_contraction2}
\diam(\P_n[x_n]) \leq C(\hat x) \lambda^{-n} \qquad \mbox{ for all } n \geq 0.
\end{equation}
\end{lemma}

\begin{proof} 
The measure 
$\hat \mu$ is invariant in forward and backward time, in particular 
$\hat \mu(A) = \hat \mu(\hat f(A))$.
Let $D\hat f^{-n}$ denote the derivative
of $f^{-n}$ restricted to an inverse branch.
By the Birkhoff Ergodic Theorem,
\begin{eqnarray*}
\lim_{n\to \infty}  \frac1n
\log |D\hat f^{-n}(\hat x)| &=& \int \log |D\hat f^{-1}| d \hat \mu =
 \int - \log |Df| d \mu = - \lambda(\mu)
\end{eqnarray*}
for $\hat \mu$-a.e. $\hat x$, so
$|Df^n(x_n)| \geq \frac1{C_0} \lambda^n$ for some $C_0 = C_0(\hat x) < \infty$.
Using  (\ref{distortion}), we find that
\[
\diam(\P_n[x_n]) \leq K \frac{|\P_0[x_0]|}{|Df^n(x_n)|} \leq C (\hat x)  \lambda^{-n}
\]
as required.
\end{proof}

As a consequence of~(\ref{distortion}) for some $K>0$
\begin{equation}\label{good_bound}
\frac{1}{K^2}\le \frac{\rho(y_{i+1},z_{i+1})\rho(\tilde{y_0},\tilde{z_0})}
{\rho(\tilde{y}_{i+1},\tilde{z}_{i+1})\rho(y_0,z_0)}\le K^2
\end{equation}
for all $y_n, \tilde y_n, z_n, \tilde z_n \in \P_n[x_n]$.
In dimension one it is sufficient to bound 
\[
\slims_{i=0}^{n-1}\|\Ad(\phi_i(y_i))\| \ d(\phi(y_{i+1}),\phi(z_{i+1})) 
\]
by
\[
\slims_{i=0}^{n-1}C \cdot K(\hat x) (\mu_u  )^{i+1} 
\lambda^{-(i+1)\alpha}\rho(y_0,z_0)^{\alpha} 
\]
in (\ref{eq_Holder}) of the proof of the main theorem. Using the two
observations  (\ref{expo_contraction2}) and (\ref{good_bound}) the
proof goes through as in Theorem~\ref{main}.

\subsection{Smooth measures}

The assumption that the Jacobian $J_{\mu}$ is H\"older is
used by e.g. Young \cite{Young1}, and enables us to apply the technique to
Gibbs measures and equilibrium states of suitable potentials.
But in the case that $\mu$ is absolutely continuous with respect to Lebesgue 
measure, it suffices to assume that the density $h \in L^1(\leb)$.
Indeed, equation (\ref{Hmc}) can be derived as follows. Let $\eps>0$
and $\eta>0$.
Choose $n$ sufficiently large that 
 $\leb(\P_n[x_n] \cap \{y\in \P_n[x_n]: d(\psi(x_n),\psi(y))<\eta\}^c)
\leq \eta \, \leb(\P_n[x_n])$. 
Boundedness of distortion gives that
\[
\frac{\leb(f^n(\P_n[x_n]\cap \{y\in\P_n[x_n]:
  d(\psi(x_n),\psi(y))<\eta\}^c))}{\leb(f^n(\P_n[x_n]))}
\]
\[ \qquad\qquad
\leq \O(1)\ \frac{\leb(\P_n [x_n] \cap   \{y\in \P_n[x_n]: 
d(\psi(x_n),\psi(y))<\eta\}^c)}{\leb(\P_n [x_n])}.
\]
Now write
\begin{eqnarray*}
\mu(\Lambda) &=&\mu(f^n (\P_n[x_n] \cap \{y\in \P_n[x_n]:
d(\psi(x_n), \psi(y))<\eta\})\\
&& + \ \mu(f^n (\P_n[x_n] \cap \{y\in \P_n[x_n]: 
d(\psi(x_n),\psi(y))<\eta\}^c))\\
&=&\int_{f^n(\P_n[x_n]\cap \{y\in \P_n[x_n]:
  d(\psi(x_n),\psi(y))<\eta\}} h(x)dx\\
&&+ \ \int_{f^n (\P_n[x_n] \cap \{y\in \P_n [x_n]: d(\psi(x_n),\psi(y))<\eta\}^c)} h(x)dx
\end{eqnarray*}
and note that since the density $h(x)\in L^1(\leb)$, and 
$\leb(f^n (\P_n[x_n] \cap \{y\in \P_n[x_n]: d(\psi(x_n),\psi(y))<\eta\}^c)) 
\leq \O(\eta) \leb(\Lambda)$ we can 
assume $\eta>0$  is such that  
$\mu(f^n(\P_n[x_n] \cap \{y\in \P_n[x_n]: 
d(\psi(x_n),\psi(y))<\eps\}^c))\leq \eps \mu(\Lambda)$.
Here we have used the  fact that for a $L^{1}(\leb)$ function $g$: 
given $\delta_1>0$, there exists a $\delta_2>0$ so that $\leb(A)<\delta_2$ implies
$\int_{A}g<\delta_1$.

\section{Refinements and Applications}

\subsection{Young towers}

We will show that Theorem~\ref{main} implies H\"{o}lder regularity for measurable
solutions to a broad class of 
cohomological equations  on   Young towers that arise in applications. 

Suppose $T:X\to X$ is a $C^{1+\gamma}$ 
mapping of a Riemannian  manifold $X$ and $\leb$ denotes   Lebesgue measure. 
Let $\rho_X$ denote the corresponding metric.
A Young tower  for $T$ has the  properties:
\begin{itemize}
\item There exists a set $\Lambda \subset X$, decomposed as
$\Lambda = \cup_j \Lambda_j \mod \leb$.
\item For each $j$, there exists $R_j \geq 1$ such that
$T^{R_j}:\Lambda_j \to \Lambda$ is bijective.
Denote the {\em induced map}  $T^{R_j}|\Lambda_j$ by $F$.
\item The distortion is bounded, i.e. there exists $K < \infty$ such that
for all $n\ge 0$
\[
\left| \frac{\Jac{DF^n_x}}{ \Jac{DF^n_y}} \right| \leq K 
\]
for all $x,y \in U$ and sets $U$ on which $F^n$ is a diffeomorphism.
\item There exists $\lambda_0 > 1$ such that 
$\min_{|v|=1}|DF_xv|/|v| \geq \lambda_0$
for  $x\in \Lambda$. 
\end{itemize}

The Folklore Theorem 
\cite{Mane} 
states that $F$ has an invariant probability $\mu$, which 
is equivalent to Lebesgue and the Radon-Nikodym derivative $h$ is bounded 
and bounded away from $0$.

The measure $\mu$ can be pulled back
to a $T$-invariant  measure $\nu$: 
\begin{equation}\label{pulled-back}
\nu(A) = \sum_j \sum_{i=0}^{R_j-1} \mu(T^{-i}A \cap \Lambda_j).
\end{equation}
The measure $\nu$ is finite if and only if
\begin{equation}\label{summable}
\mathcal{R} := \int_{\Lambda} R\  d\mu = \sum_j R_j \mu(\Lambda_j) < \infty.
\end{equation}

This set-up can be viewed as a Markov system as follows:
\begin{itemize}
\item $M$ is the disjoint union 
$\sqcup_j \sqcup_{i=0}^{R_j-1} (\Lambda_j,i)$
where each $(\Lambda_j,i)$ is a copy of $\Lambda_j$.
The set $\Lambda=\cup_j (\Lambda_j,0)$  is called the {\em base} of 
the tower. Each set $(\Lambda_j,i)$ is a component $M_k$ of $M$.
\item For the metric $\rho$ on $M$, there are at least two choices.
Take $(x,i) \in (\Lambda_j,i)$ and 
$(\tilde x,\tilde i) \in (\Lambda_{\tilde j},\tilde i)$
\begin{itemize}
\item[$\diamond$]
$\rho_1((x,i), (\tilde x, \tilde i)) = 
\left\{ \begin{array}{ll}
\rho_X(T^i(x), T^{\tilde i}(\tilde x)) & \mbox{ if }  
j = \tilde j \mbox{ and } i = \tilde i;\\
1 & \mbox{ otherwise.} \\
\end{array} \right.
$\\[0.2cm]
This metric is induced from the metric $\rho_X$ on $X$. 
The metric $\rho_1$ is used in Corollary~\ref{original_system}.\\
\item[$\diamond$] 
$
\rho_2((x,i), (\tilde x, \tilde i)) = 
\left\{ \begin{array}{ll}
\rho_X(x, \tilde x) & \mbox{ if }  j = \tilde j \mbox{ and } i = \tilde i;\\
1 & \mbox{ otherwise.} \\
\end{array} \right.
$\\[0.2cm]
This metric is the {\em tower metric}
and is induced from the metric on the base $\Lambda$.
\end{itemize}
\item
Define $f:M \to M$ as 
\[
f(x,i) = \left\{ \begin{array}{ll}
(x,i+1) & \mbox{ if } x \in \Lambda_j \mbox{ and } i < R_j-1,\\
(T^{R_j}x,0) & \mbox{ if } x \in \Lambda_j \mbox{ and } i = R_j-1.
\end{array} \right.
\]
There is a projection $\pi:M \to X$ given by $\pi(x,i) = T^i(x)$
and $\pi \circ f = T \circ \pi$.
\item Extend the definition of $\mu$ from $\Lambda$ to $M$
as $\mu((A,i)) = \mu(A)$ for each measurable set $A \subset \Lambda_j$
and $0 \leq i \leq R_j-1$.
By \eqref{summable} $\mu$ is finite and hence can be normalized.
\end{itemize}

Let $\P_0$ denote the partition of $M$ into sets $(\Lambda_j,i)$ and 
set $\P_n=\bigvee_{i=0}^{n-1}f^{-i}(\P_0)$.  For each $x\in M$
let $\P_n[x]$ be the partition element (cylinder set) in $\P_n$
containing $x$.
Let $(\hat M, \hat f, \hat \mu)$ denote the natural extension 
of $(M,f,\mu)$. For $\hat x=(x_0,x_1,\ldots, x_n,\ldots)$ let
$n_0 < n_1 < \dots$ denote the indices such that
$x_{n_i}$ belongs to the base.

\begin{lemma}\label{hyp_remove}
Let $(M,f,\mu)$ be a Young tower satisfying \eqref{summable}.
Assume one of the following three conditions:
\begin{itemize}
\item  $\rho = \rho_1$ and there is $\delta_0 > 0$ such that
\begin{equation}\label{bounded_contraction}
\| DT^{R_j - k}(T^k(x))\| \geq \delta_0 \mbox{ for all } j, x 
\in \Lambda_j \mbox{ and } 0 \leq k < R_j.
\end{equation}
\item  $\rho = \rho_2$, or \\
\item  $\log^+\|Df^{-1}\|\in L^1 (\mu)$, where the derivative is taken
with respect to the metric used.
\end{itemize}
Then for $\hat \mu$-a.e. $\hat x$ 
there exists a number $K(\hat x)$ and $\lambda>1$  such that 
\[
\rho(y_n,z_n)<K(\hat x) \lambda^{-n} \rho(y_0,z_0)
\]
for all $n\ge 0$ and $y_n,z_n\in \P_n [x_n]$.
\end{lemma}
 
\begin{proof}
Write $R(x) = R_j$ for $x \in \Lambda_j$, so by 
Birkhoff's and Kac's Theorems
${\mathcal R}^{-1} = \lim_{k \to \infty} \frac1k \sum_{i = 1}^k R(x_{n_i})$
for $\hat\mu$-a.e. $\hat x$.
Let $\tilde \Lambda_i = \cup \{ \Lambda_j \ : \ R_j = i\}$,
so \eqref{summable} gives $\sum_i i \mu(\tilde \Lambda_i) < \infty$.
Given $\eps > 0$, we have
\begin{eqnarray*}
\sum_k \hat \mu\{ \hat x \ : \ R(x_{n_{k+1}}) > \eps n_k \} &\leq&
\sum_k \hat \mu\{ \hat x \ : \ x_{n_{k+1}} \in \cup_{i \geq \eps k} 
\tilde \Lambda_i \} \\
 &\leq&\sum_k \sum_{i \geq \eps k} \mu(\tilde \Lambda_i) \\
 &\leq& \lceil \frac1{\eps} \rceil \ \sum_k k \mu(\tilde \Lambda_k) < \infty.
\end{eqnarray*}
The Borel-Cantelli Lemma therefore implies that for $\hat \mu$-a.e. $\hat x$,
there is $k_0 = k_0(\hat x)$ such that $n_{k+1}-n_0 \leq (1+\eps) (n_k-n_0)$ 
for all $k \geq k_0$. Take $n_k \leq n < n_{k+1}$, then 
\[
\| Df^{-(n-n_k)}(x_{n_k}) \| \leq \left\{
\begin{array}{ll} 
\delta_0^{-1} & \mbox{ if } \rho = \rho_1 
\mbox{ by \eqref{bounded_contraction};} \\
\lambda_0^{-1} < 1 & \mbox{ if } \rho = \rho_2.
\end{array} \right.
\]
Write $B(\hat x) = \| Df^{-n_0}(x_0) \|$.
Direct calculation gives
\begin{eqnarray}\label{lya_tower}
\lim_{n \to \infty}
\left(\|\prod_{i=0}^{n-1}Df^{-1} (x_i)\|\right)^{\frac{1}{n}} 
&\leq&\lim_{k\to\infty} \left( \frac{B(\hat x)}{\delta_0} 
\|\prod_{i=0}^{n-1}Df^{-1} (x_i)\|\right)^{\frac1{(1+\eps)n_k}} \notag\\
&=&\lim_{k \to \infty} \left( \frac{B(\hat x)}{\delta_0} \|\prod_{i=1}^k
DF^{-1} (x_{n_i})\|\right)^{\frac{1}{n_0 + (1+\eps)\sum_{i=1}^{k} R(x_{n_i})} }
\notag \\
&=&\lim_{k \to \infty} \left( \bigl(\frac{B(\hat x)}{\delta_0} \|\prod_{i=1}^{k}
DF^{-1} (x_{n_i})\|\bigr)^\frac{1}{k}\right)^{
\frac{k}{n_0 + (1+\eps)\sum_{i=1}^k R(x_{n_i})} } \notag \\
&\le & \lambda_0^{-\frac{1}{(1+\eps)\mathcal R}}.
\end{eqnarray}
Because $\eps > 0$ was arbitrary, we get 
$\lim_{n \to \infty}
\left(\|\prod_{i=0}^{n-1}Df^{-1} (x_i)\|\right)^{\frac{1}{n}}
\leq \lambda_0^{-\frac{1}{\mathcal R}}$ as well.

Alternatively, assume that $\log^+\|Df^{-1}\|\in L^1(\mu)$. 
The estimate of \eqref{lya_tower} holds for the sequence $(n_k)$. 
Oseledec's Multiplicative Theorem then implies that
$\lim_{n \to \infty}
\left(\|\prod_{i=0}^{n-1}Df^{-1} (x_i)\|\right)^{\frac{1}{n}}$
exists $\hat \mu$-a.e. $\hat x$. This limit equals the limit along the 
subsequence, so
$$\lim_{n \to \infty}
\left(\|\prod_{i=0}^{n-1}Df^{-1} (x_i)\|\right)^{\frac{1}{n}} \leq 
\lambda_0^{\frac{1}{\mathcal R}}.$$

Since these estimates hold uniformly over $y_i,z_i\in \P_i [x_i]$, 
we obtain $\rho(y_n,z_n)<K(\hat x) \lambda^{-n} 
\rho(y_0,z_0)$ for $\lambda=\lambda_0^{\frac{1}{\mathcal R}}$.
\end{proof}

\begin{remark} The condition $\log^+\|Df^{-1}\|\in L^1(\mu)$ relates to 
the reference measure $m$ on $X$ as follows: if $\mu$ has density $h$, 
and $h \in L^{1+\eps}(m)$ for some $\eps > 0$, then by the H\"older
inequality, $\log^+\|Df^{-1}\|\in L^q(m)$ for $q = \frac{1+\eps}{\eps}$
implies that $\log^+\|Df^{-1}\|\in L^1(\mu)$.
\end{remark}

\begin{theorem}\label{young_tower}
Let $(M,f,\mu)$ be a Young tower with base $\Lambda$,
satisfying \eqref{summable}.
Assume that $\log^+\|Df^{-1}\|\in L^1(\mu)$ or that 
\eqref{bounded_contraction} holds.
Let  $\phi:M \to G$ be  a Lie group valued $\alpha$-H\"older
observable and define  
\[
\mu_u := 
\lim_{n\to\infty}\left( 
\sup_{(x,0)\in M} \| \Ad (\phi(f^{n-1} (x,0)\ldots
 \phi(f(x,0))\phi ((x,0))) \| \right)^{\frac{1}{n}}.
\]
If $1 \leq \mu_u < \lambda^{\frac{\alpha}{\mathcal{R}}}_0$, then
any $\mu$-measurable solution to the cohomological equation
\[
\phi = (\psi \circ f) \psi^{-1} \qquad \mu\mbox{-a.e.}
\]
 has a version which is  $\alpha$-H\"older on $\Lambda$.
\end{theorem}

\begin{remark}
If $G$ is Abelian, compact or nilpotent  then $\mu_u=1$ and the spectral condition 
 $1 \leq \mu_u < \lambda^{\frac{\alpha}{\mathcal{R}}}_0$
is automatically satisfied.
\end{remark}

\begin{remark}
For any $n$ there exists  an $\alpha$-H\"{o}lder version on 
$\cup_{j=0}^n f^j (\Lambda)$ but the H\"older constant may
increase with $n$.

\end{remark}

\begin{remark}
The same result holds for  solutions $\psi$ to equation~\eqref{eq-cohomology1},
equation~\eqref{eq-cohomology2} or  equation~\eqref{eq-cohomology3}.
\end{remark}

\begin{proof}
First observe that $\mu$ is an ergodic $f$-invariant measure.
If $\rho = \rho_2$, and $x\in \Lambda_i$,
then $Df (f^i (x))= \mbox{Id}$ if $i<R_i-1$ 
and $Df(f^{R_i(x)})=DF(x)$.
If $\rho = \rho_1$, then $Df(f^i(x)) = DT(T^i(x))$.
In either case, $Df(f^{R_i(x)})=DF(x)$.

A computation similar to \eqref{lya_tower} yields
that
for $\mu$-a.e. $x \in \Lambda$:
\begin{equation}\label{eq:lya} 
\lambda(\mu) :=\lim_{n\to\infty}
\left(\|\prod_{i=0}^{n-1}Df (f^i (x,0))\|\right)^{\frac{1}{n}} 
= \lambda_0^\frac{1}{\mathcal R}.
\end{equation}
Therefore the assumption $\mu_u < \lambda^{\frac{\alpha}{\mathcal{R}}}_0$ 
implies partial hyperbolicity (\ref{PH}).
Lemma~\ref{hyp_remove} shows that
condition~(\ref{expo_contraction}) holds. Hence 
the proof of Theorem~\ref{main} applies and
we obtain a H\"{o}lder version on each $(\Lambda_j, i)$. 
Using the argument from 
Corollary~\ref{cor_extension}
we obtain Liv\u{s}ic regularity on the base $\Lambda$ as well.
\end{proof}

\begin{corollary}\label{original_system}
Suppose $(T,X)$ is modelled by a Young tower $(M,f, \mu)$ over base 
$\Lambda \subset X$, satisfying \eqref{summable}.
Assume that $\log^+\|Df^{-1}\|\in L^1(\mu)$ or that 
\eqref{bounded_contraction} holds.
Let $\nu$ be the $T$ invariant ergodic  pulled back measure,
as in equation \eqref{pulled-back}.
Let $\tilde \phi:X \to G$ be $\alpha$-H\"older and satisfy
$\mu_u < \lambda_0^{\frac{\alpha}{\mathcal R}}$. 
Then for any $0 \leq k \leq {\mathcal R}_j$,
 any $\nu$-measurable solution to the cohomological equation
$\tilde \phi = (\tilde \psi \circ T) \tilde \psi^{-1}$ has a
version which is H\"older on $T^k(\Lambda_j)$.
\end{corollary}

\begin{proof}
Suppose   $\tilde{\phi}:X\to G$
is H\"{o}lder and $\tilde{\phi}=(\tilde{\psi}\circ T) \tilde{\psi}^{-1}$,
$\nu$-a.e.
Let $\pi:M \to X$ be defined as $\pi((x,i)) = T^i(x)$,
so that $T \circ \pi = \pi \circ f$.
Then $\tilde \phi, \tilde \psi:X \to G$ lift
to the tower as $\phi = \tilde \phi \circ \pi$, $\psi=\tilde \psi \circ \pi$
and satisfy $\phi=(\psi\circ f){\psi}^{-1}$,
$\mu$-a.e. Moreover, $\phi$ is $\alpha$-H\"older
with respect to the metric $\rho_1$. Since
derivatives $Df( f^i(x) )$ on the tower agree with derivatives
$DT(T^i(x))$ on $(X,T)$, 
\[
\left(\|\prod_{i=0}^{n-1}Df (f^i(x))\|\right)^{\frac{1}{n}}\ge
\lambda_0^{\frac{\alpha}{\mathcal R}}
> \mu_u
\]
for $\mu$-a.e. $x$, verifying (PH) on $(M,f,\mu)$.

Now Theorem~\ref{young_tower} gives 
$\alpha$-H\"older (with respect to $\rho_1$)
solution $\psi'$ to the cohomological equation
$\phi=(\psi '\circ T){\psi '}^{-1}$,
such that $\psi'=\psi$ $\mu$-a.e.
Since $\psi'$ takes the same value on every point of $\pi^{-1}(x)$
for $\nu$-a.e. $x$
the projection $\tilde \psi' = \psi \circ \pi$ is well-defined,
$\alpha$-H\"older and satisfies 
$\tilde{\phi}=(\tilde{\psi '}\circ T) \tilde{\psi '}^{-1}$,
$\nu$-a.e.
\end{proof}

\begin{remark}
A priori, a larger class
of observables $\phi:M \to G$ is H\"older with respect to the 
tower metric $\rho_2$ than with respect to $\rho_1$.
For metric $\rho_2$, however, the
projection $\pi:M \to X$, $(x,i) \mapsto T^i(x)$ need not
preserve the H\"olderness of solutions
of the cohomologous equation on $(M,f)$.
If the projection $\pi$ preserves continuity or is only used on the 
base $\Lambda$, this may still suffice for applications.
\end{remark}

\subsection{The Manneville-Pomeau family}\label{manneville_pomeau}

In the previous results, the assumptions that $\mu$ is finite
and/or has positive Lyapunov exponents can be weakened
for some group extensions. 

For the Manneville-Pomeau family the Jacobian $J_{\mu}$ of $\mu$ is H\"older:
i.e. there exists $C$ and $\gamma \in (0,1]$ such that
\begin{equation}
\left| \frac{J_\mu(x)}{J_{\mu}(y)} - 1 \right| \leq 
C \cdot \rho(f(x),f(y))^{\gamma}.
\end{equation}

To prove Lemma~\ref{dist_Jac}, assumption 
\eqref{expo_contraction} can be weakened to
\begin{equation}\label{eq_Pn_summable}
\sum_{n \geq 0} \rho(y_n,z_n)^{\gamma} < K(\hat x) \rho(y_0,z_0)^{\gamma},
\end{equation}
where $\gamma > 0$ is the H\"older exponent in (\ref{Jac_Holder}).
In this case $\lambda(\mu)$ can be $0$, and (\ref{PH}) fails,
but if $G$ is Abelian or compact, or if $G$ is nilpotent
with a dominated growth rate so that for inverse branches
\begin{equation}\label{eq:nilpotent}
\sum_{n\geq 0}\|\Ad(\phi_n(y_0))\| \ 
d(\phi(y_{n+1}),\phi (x_{n+1}))^{\alpha}<\infty,
\end{equation}
then the estimates (\ref{eq_Holder}) hold.
Furthermore, the pulled back measure $\nu$ can be at most $\sigma$-finite
if (\ref{summable}) fails.
According to equation~\eqref{eq:lya},  
this implies that $\lambda(\mu) \leq 1$.

This scenario is found in  the well-known Manneville-Pomeau maps.
This is a family of maps on $[0,1]$ which have a neutral fixed point
(where we take $\rho$ to be Euclidean distance), 
parameterized by parameter $p\in(0,\infty).$ 
For $p\geq 1$ these maps admit a $\sigma$-finite absolutely continuous 
invariant measure. However
the measure is not a probability measure, and the map has zero 
Lyapunov exponents for Lebesgue
almost all initial conditions. We have the  following regularity result in this setting:
\begin{theorem}\label{indifferent}
Consider the Manneville-Pomeau map (for $p \geq 0$):
\[
T(x) = \left\{ \begin{array}{ll}
x + 2^px^{1+p} & \mbox{ if } x \in [0,\frac12), \\
2x-1 & \mbox{ if } x \in [\frac12, 1].
\end{array} \right.
\]
If $G$ is an Abelian or compact group and
$\tilde \phi:[0,1] \to G$ is a $G$-valued 
$\alpha$-H\"older observable
for $\alpha > \frac{p}{1+p}$, then 
any Lebesgue measurable solution to the cohomological equation 
$\phi = (\psi \circ T) \psi^{-1}$ has  a  version which is 
$\alpha$-H\"older on $[0,1]$.
\end{theorem}

\begin{remark}
The same regularity is forced upon solutions $\psi$ to equation~\eqref{eq-cohomology1},
equation~\eqref{eq-cohomology2} or  equation~\eqref{eq-cohomology3}.
\end{remark}

\begin{proof}
For $p \in [0,1)$, $T$ has an invariant probability $\nu \ll \leb$,
whereas $\nu$ is only $\sigma$-finite if $p \geq 1$.
However, the first return map $F$ to $[\frac12, 1]$ has always an invariant
probability $\mu \ll \leb$, and $F^n$ has bounded distortion, independently
of $n$.
Therefore we can use Theorem~\ref{young_tower} 
provided we can show that
$\phi|(\Lambda_j,i) = \tilde \phi \circ \pi^{-1}|\Lambda_j = 
\phi \circ T^i|\Lambda_j$ is H\"older.
This is done as follows.
It is not hard to check that 
$T^{-n}(\frac12) = \frac12 (pn)^{-\frac1p} + o(n^{-\frac1p})$,
where $T^{-n}$ indicates the $n$-th inverse of the left branch of 
$T$. Consequently, $\diam(T^{-n}(J)) = \O(n^{-\frac{1+p}{p}})$.

Even if  $\nu$ is $\sigma$-finite, $\nu([\frac12, 1]) < \infty$.
Therefore the subset 
$\hat M_0 := \{ \hat x \in \hat M \ : \ x_0 \in [\frac12,1]\}$
has finite $\hat \mu$-measure.
Suppose that $\hat x \in \hat M_0$ 
is a backward orbit 
(chosen as in Theorem~\ref{main}), 
and $\hat y, \hat z \in M$ are such that $y_n, z_n \in \P_n[x_n]$ for 
each $n$. Let $n_0 = 0$ 
and $n_k = \min\{ n > n_{k-1} \ | \ x_n\in \Lambda = [\frac12,1]\}$.
Since $|DF| \geq 2$, $\rho(y_{n_k}, z_{n_k}) \leq 2^{-k}$, and hence,
using the fact that $\alpha > \frac{p}{1+p}$,
\begin{eqnarray*}
\sum_{n \geq 0} d(\phi(y_n) ,\phi(z_n)) &\leq&
\sum_{k \geq 0} \sum_{n=n_k}^{n_{k+1}-1} \O(1) \rho(y_n , z_n)^{\alpha} \\
&\leq&\sum_{k \geq 0} \rho(y_0, z_0)^{\alpha} 2^{-\alpha k}
\sum_{n=0}^{n_{k+1}-n_k-1} \!\! \O(1) \diam( T^{-n}([\frac12, 1]) )^{\alpha} \\
&\leq& \rho(y_0,z_0)^{\alpha} \sum_{k \geq 0} 2^{-\alpha k} 
\sum_{n=0}^{n_{k+1}-n_k-1} \!\! \O(1) n^{-\alpha \frac{1+p}{p}} \\
&\leq& \O(1)  \rho(y_0,z_0)^{\alpha}.
\end{eqnarray*}
This calculation replaces (\ref{eq_Holder}) in the proof of 
Theorem~\ref{main}. Continuing the proof as in Theorem~\ref{main}, we get that
$\psi$ has an $\alpha$-H\"older version on $[\frac12, 1]$.
Because $T:[\frac12,1] \to [0,1]$ is smooth by using the 
cohomological equation
$\psi \circ T(x) = \phi(x) \cdot \psi(x)$, it follows that $\psi$ has an 
$\alpha$-H\"older version on $[0, 1]$, cf. Corollary~\ref{cor_extension}.
\end{proof}

\begin{remark} The arguments above can also be used to establish 
H\"older regularity for nilpotent groups $G$,
provided (\ref{eq:nilpotent}) holds. 
Having $\mu_u=1$ alone is not enough to establish Liv\u{s}ic regularity, but
if $\|\Ad(\phi_n(y_0))\|$ grows at a  rate `dominated' by  the polynomial contraction of 
$d(\phi(y_{n+1}),\phi (x_{n+1}))^{\alpha}$, which in our case is 
$\O(n^{-\frac{\alpha(1+p)}{p}}),$ then the conclusion of Theorem~\ref{indifferent} will
remain valid.
\end{remark}

\subsection{Interval Maps with Critical Points}
So far, the Markov systems used in the examples were Young towers,
even though we allowed a $\sigma$-finite measure for the 
Manneville-Pomeau map.
Theorem~\ref{main} also applies to different kinds of Markov systems.
In this subsection we discuss the consequences of the theory
to systems modelled by so-called Hofbauer towers.
The metric $\rho$ is Euclidean distance in this section.

An interval map $T:[0,1] \to [0,1]$ is called piecewise 
continuous (piecewise $C^r$)
if there exists a finite set of points
$0 = a_0 < a_1 < \dots < a_{k-1} < a_k = 1$ such that
$T|(a_{i-1}, a_i)$ has a continuous ($C^r$) extension to $[a_{i-1}, a_i]$. 
We call  $a_0 , \dots, a_k$ the {\em critical} points
of $T$ and denote this set by $\Crit$.
Let $\P_1$ denote the partition $\{ [a_0,a_1], \dots, [a_{k-1},a_k]\}$,
and $\P_n = \bigvee_{i=0}^{n-1} T^{-i}(\P_1)$
be the partition into $n$-cylinders. 
For $x \in I \setminus \cup_{i=0}^{n-1} T^{-i}(\Crit)$, let $\P_n[x]$ 
denote the $n$-cylinder containing $x$.

Due to the presence of critical points, $T^{-1}$ need not be 
differentiable at boundary points of $T(P)$, $P \in \P_1$, and hence 
boundedness of  distortion 
(i.e. condition \eqref{distortion}) cannot be realized globally.
For this reason, we
call $J$ a {\em core} interval if it is compactly 
contained in $T^n(P)$ for some $n \geq 1$ and $P \in \P_n$.
For example, if $T (x) = 1-ax^2$ is a non-renormalizable unimodal map,
then any interval compactly contained in $[T^2(0), T(0)]$
is a core interval.

\begin{theorem}\label{pw}
Let $T$ be a piecewise $C^3$ interval map onto the unit interval
with negative Schwarzian derivative.
Let $\nu$ be a $T$-invariant probability measure 
such that $\int \log |T'| \ d\nu > 0$. Assume that the Jacobian $J_{\nu}$
is H\"older or $\nu \ll \leb$.
Let $\phi$ be a Lie group valued
piecewise $\alpha$-H\"older observable, with discontinuities
(if any) only at the points $a_i$, and satisfying 
the partial hyperbolicity condition (\ref{PH}).
If $\phi = (\psi \circ T) \psi^{-1}$ for some $\nu$-measurable function,
then for every {\em core} interval $J$ with $\nu(J) > 0$,
$\psi|J$ has an $\alpha$-H\"older version.
\end{theorem}

\begin{remark}
Assume that the core interval $J$ is compactly contained in $T^n(P)$
for some $P \in \P_n$, and $n \geq 1$.
Although the H\"older exponent is independent of $J$,
the H\"{o}lder coefficient in the H\"older property of the H\"older version 
of $\psi|J$ depends on
$\min\{ \frac{\tiny \diam(J)}{\tiny \diam(L)}, 
\frac{\tiny \diam(J)}{\tiny \diam(R)} \}$, 
where $L$ and $R$ are 
the components of $T^n(P) \setminus J$.
As this minimum tends to $0$, the H\"{o}lder coefficient tends to infinity.
\end{remark}

\begin{remark}
The assumption that $T$ is onto is not  a severe restriction.
Since $\nu$ is assumed to have a positive Lyapunov exponent,
$\nu$ cannot be supported on a non-repelling periodic orbit,
so $\cap_j T^j([0,1])$ is a non-trivial interval.
For the same reason, $\nu$ cannot be supported on a wandering 
interval, i.e. an interval $W$ such that $T^n(W) \cap T^m(W) = \emptyset$
for all $n \neq m \geq 0$.
By restricting and rescaling, we can assume that $T:[0,1] \to [0,1]$ is onto.
\end{remark}

\begin{remark}
The same result holds for solutions $\psi$ to equation~\eqref{eq-cohomology1},
equation~\eqref{eq-cohomology2} or  equation~\eqref{eq-cohomology3}.
\end{remark}

\begin{proof}
We  construct a Markov system, introduced by
Hofbauer~\cite{Hofbauer} as the {\em canonical Markov extension}
and sometimes  called  a {\em Hofbauer tower}.
This is the system $(M,f)$, where $M$ is a disjoint union of 
closed intervals. 
We call $B = [0,1]$ the {\em base of the tower}. Then
\[
M = B \sqcup 
\left( \sqcup_{n \geq 1} \sqcup_{P \in \P_n} \overline{T^n(P)} \right)/\sim
\]
where the sets $T^n(P) \sim T^{n'}(P')$ if they are the same interval.
Let $\pi:M \to [0,1]$ be the natural projection.
The action $f$ is defined on $M$ as follows. If $x \in M$ belongs to 
the component $D$, then
\[
f(x) = \pi^{-1}(T(\pi(x))) \cap \tilde D,
\]
where component $\tilde D := \overline{T(D \cap \P_1[\pi(x)])}$ 
is again a component of $M$.
Obviously $(M,f)$ is Markov and $T \circ \pi = \pi \circ f$.
Due to the Markov property, the following holds for any component $D$ of $M$:
\[
T^n(\P_n[x]) = \pi(D) 
\mbox{ if and only if } f^n(\pi^{-1}(x) \cap B) \in D.
\]
If $\nu$ is $T$-invariant, then we can construct a measure $\mu$ as follows.
Let $\mu_0$ be the measure $\nu$ lifted to the level $B$
and set $\mu_n = \frac{1}{n+1} \sum_{i=0}^{n} \mu_0 \circ f^{-i}$.
Clearly $\nu = \mu_n \circ \pi^{-1}$ for each $n$.
As shown in \cite{keller}, $\mu_n$ converges vaguely. 
We call the limit measure $\mu$.
If $\nu$ is ergodic, then $\mu$ is either a probability measure on $M$,
in which case we call $\mu$ {\em liftable},
or it is identically $0$ on $M$.
In this case the mass ``has escaped to infinity''.
Keller's result \cite[Theorem 3]{keller} states 
if $\int \log |T'| \ d\nu > 0$, then $\nu$ is 
liftable to an invariant measure $\mu$ on the Markov extension. 
Moreover, $\mu$ is ergodic if $\nu$ is.

If $J$ is a core interval and $\nu(J) > 0$, then there is some
level $D \in M$ compactly containing a lifted copy $\tilde J := D \cap \pi^{-1}(J)$,
and $\mu(\tilde J) > 0$.
Let $\delta > 0$ be such that $D$ contains a $\delta |J|$-neighbourhood of $\tilde J$.
Using negative Schwarzian derivative and the Koebe principle, see \cite{dMvS},
we find that for every $x \in M$ with $f^n(x) \in \tilde J$,
that $f^n$ has bounded distortion on the component of $f^{-n}(\tilde J)$ 
containing $x$. In fact, the distortion depends only on $\delta$.

Finally, we will show that $\psi|J$ has a H\"older version.
 $\tilde \phi = \phi \circ \pi$ is an observable on the Markov extension.
The coboundary $\psi$ lifts to a coboundary $\tilde \psi = \psi \circ \pi$.
Therefore we can apply Theorem~\ref{main} to it to find a version of
$\tilde \psi$ that is $\alpha$-H\"older. Projecting it back to the interval, 
we find the desired $\alpha$-H\"older version of $\psi$.
(Note that since $\tilde \psi$ takes the same value on every point in 
$\pi^{-1}(x)$, we find that the H\"older version of $\psi$ 
does not depend on the level $D \in M$ that we lift $J$ to.)
\end{proof}

Let $T:I \to I$ be a $C^3$ S-multimodal map having 
critical set $\Crit$. Assume each critical point $c \in \Crit$ has order
$\ell_c$, $1 < \ell_c < \infty$. 
We assume for simplicity that $T$ is locally eventually
onto, i.e., for every 
non-degenerate subinterval $U \subset I$, $T^n(U) = I$
for some $n \geq 0$. This excludes that $g$ is renormalizable, or has a 
non-expanding periodic orbit.
(Also wandering intervals are excluded, but this is already
a corollary of the smoothness, see \cite{dMvS}.)
In this case the $T$-invariant measure $\nu$ that we will be considering
is supported on $I$.

\begin{theorem}\label{Hofb}
Assume $T:I \to I$ is $C^3$ multimodal with negative Schwarzian derivative 
and non-flat critical points. Write $c_n = T^n(c)$ and 
$\ell_{\max} = \max\{ \ell_c \ : \ c \in \Crit\}$.
Assume that $T$ satisfies the summability condition
\begin{equation}\label{sum}
\sum_{c \in \Crit} \sum_{k \geq 1} |DT^{k-1}(c_1)|^{-1/\ell_{\max}} < \infty,
\end{equation}
and hence possesses an acip $\nu$ (cf. \cite{NvS,BvS}).
Let $\phi$ be a piecewise $\alpha$-H\"older $L^1(\nu)$ observable,
with discontinuities of types (1)-(3) of Subsection~\ref{sec:sing}, 
at critical points only. 
If $\phi = \psi \circ T - \psi$ for some $\nu$-measurable function,
then for every core interval $J$,
$\psi|J$ is $\tilde \alpha$-H\"older for any $\tilde \alpha \in (0,\alpha)$.
\end{theorem}

\begin{proof}
Assume as above that $T$ is locally eventually onto.
We 
build a Markov extension $(M,f)$ as in 
Theorem~\ref{pw}. Since $T$ admits an acip $\nu$ 
(with necessarily positive Lyapunov
exponent, cf. \cite{Keller_Hopf} and also \cite{BL}), 
it can be lifted to an acip $\mu$ on $(M,f)$.
Let $D$ be any level in $M$ that compactly contains a lifted copy 
$\tilde J = D \cap \pi^{-1}(J)$ of $J$ and such that $\nu(J) > 0$.

Let $(\hat M, \hat f, \hat \mu)$ be the natural extension 
of $(M,f,\mu)$. By the Koebe principle, we have a uniform distortion
bound for $f^n|\P_n[x_n] \cap f^{-n}(J)$ for each $x \in \tilde J$
(i.e. $J$ is a core interval as introduced before).
 
The next thing to check is that there are sufficiently many backward
paths that avoid passing close to the singularities of $\phi$ at $\Crit$.
This argument is similar to the one in the proof of Theorem~\ref{nonholder}.
The proof that $(I,f)$ has an acip was given in the multimodal case
in \cite{BvS}, based on the well-known result of Nowicki \& van 
Strien \cite{NvS}. An important estimate in \cite{BvS} is
that for some $C = C(T) < \infty$, 
\[
\nu(A) \leq C \cdot \leb(A)^{1/\ell}
\mbox{ for all measurable }
A \subset I.
\]
It follows that if $B_n = ( c-n^{-2\ell}, c+n^{-2\ell} )$, then
\[
\nu(B_n) \leq C \cdot \leb(B_n)^{1/\ell} \leq \frac{2C}{n^2},
\]
and hence, lifted to the tower:
\[
\mu(T^{-n}(\tilde J) \cap B_n) \leq \nu(B_n) \leq \frac{2C}{n^2}.
\]
It follows from the Borel-Cantelli Lemma that for 
$\hat \mu$-a.e. $\hat x \in \hat M$,
there exists $N = N(\bar x)$ such that $x_n \notin B_n$ 
for all $n \geq N$.
From now on, the argument is the same as in Theorem~\ref{nonholder}.
\end{proof}

In view of questions raised in e.g. \cite{BV}, 
we are particularly interested in 
the potential $\phi := \int_I \log |T'| d\mu - \log |T'|$.
This potential is smooth, except for logarithmic singularities at the 
critical points.
Theorem~\ref{Hofb} shows that any $\nu$-measurable solution of the
cohomological equation $\phi = \psi - \psi \circ T$
has a version which is $\alpha$-H\"older for any $\alpha \in (0,1)$
on each interval that is compactly contained in $[c_2,c_1]$.
We can apply it to the quadratic family $f_a(x)= 1-ax^2$.
It is known that for $\leb$-a.e. $a \in [0,2]$,
$f_a$ has either an attracting periodic orbit, or
has a positive Lyapunov exponent at the critical value, and both parameter
sets have positive measure.
Below we make a weaker assumption on the growth rate of derivatives
along the critical orbit.

\begin{corollary}\label{CorPhi}
Let $f : I \to I$ be a $C^3$ S-unimodal map with critical 
order $\ell < \infty$, satisfying the summability
condition $\sum_n |Df^n(c_1)|^{-1/\ell} < \infty$ (and hence possessing an 
acip $\nu$).
Then $\phi = \log |f'| - \int \log |f'| d\nu$ can only
be a measurable coboundary if there exists a periodic interval
$J \subseteq I$ of period $k$
such that $c \in J$ and $J = [f^{2k}(c), f^k(c)]$.
\end{corollary}

\begin{proof}
First assume that $f$ is nonrenormalizable. 
In this case, $\psi$ is H\"older 
continuous and hence bounded on any interval compactly contained
in $[c_2,c_1]$. Since $\phi$ is bounded except at $c$, we easily derive 
that $\psi$ has to be unbounded at every forward image of the 
critical point. This is only possible if 
$\cup_{n \geq 1} f^n(c) \subset \{ c_1, c_2 \}$.

Next assume that $f$ is finitely renormalizable, say $J \owns c$ is a 
periodic interval,
$f^k(J) \subset J$, $f^k(\partial J) \subset \partial J$ where $k \geq 2$ is
the period of renormalization, and $f^k|J$ is unimodal and nonrenormalizable.
Then the above argument shows that $\cup_{n \geq 1} f^n(c) \cap J^{\circ} = 
\emptyset$.
Therefore, $J = [f^{2k}(c), f^k(c)]$, and $f^k|J$ is conjugate to 
$x \mapsto 1-2x^2$. 
\end{proof}

In the context of parametrized S-unimodal families $f_a$, the
conclusion of this corollary is that under the summability
condition  $\sum_n |Df^n(c_1)|^{-1/\ell} < \infty$, $\phi$ is not 
a measurable coboundary except, possibly, for countably many parameter values.
If $f_a$ is infinitely renormalizable, then the summability
condition fails. If $f_a$ is finitely renormalizable, then equivalently,
there is a smallest periodic interval $J \owns c$ with $f^k(J) \subset J$.
Only when the renormalized map $f^k:J \to J$ is conjugate to
the full unimodal map $x \mapsto 1-2x^2$ (and for each $k \geq 1$ and
configuration of $J, f(J), \dots, f^{k-1}(J)$, this usually holds for
only one parameter), 
then it is possible that $\phi$ is a measurable coboundary.
However, even in this situation, it seems extremely unlikely that
$\phi$ is a measurable coboundary, as it requires that the multipliers 
of all periodic orbits are the same.

{\bf Example}: 
If $f(x) = 1-ax^2$ is the quadratic map, 
then $\phi = \log |f'| - \int \log |f'| \ d\nu$ is a 
measurable coboundary for $a = 2$, but not for the parameter 
$a \approx 1.54368901$ at which $f$ is renormalizable of period $2$, and 
$f^3(c)$ is the orientation reversing fixed point.

\begin{proof}
For $a = 2$, then $f$ is a Chebychev polynomial, and hence
$h(x) = -\cos 2x$ conjugates $T$ with the tent map 
$T(x) = \min\{ 2x , 2(1-x) \}$:
$h \circ T = f \circ h$.
It follows that $\phi = \log |f'| - \log 2 = \psi - \psi \circ f$ for
$\psi = -\log |h' \circ h^{-1}|$. 

Next assume that $f$ has a $k$-periodic interval $J$ 
as in the proof of Corollary~\ref{CorPhi}, such that
$F := f^k|J$ is conjugate to $x \mapsto 1-2x^2$.
Applying the cohomological equation to $F$, we get
$\log|F'| - k \int \log |f'| \ d\nu = \psi - \psi \circ F$.
If $h:[0,1] \to J$ is defined by $\psi = -\log |h' \circ h^{-1}|$, then
$F \circ h = h \circ T$ for the tent map $T$ as above. In particular,
$k \int \log |f'| \ d\mu = \log 2$.
By the cohomological equation, each periodic point 
$y \in J^{\circ}$ must have multiplier $\log 2$.

For $k = 2$, $F|J$ has an orientation reversing fixed point $q_1$, and
$\{ q_1, q_2 := f(q_1)\}$ is the corresponding period $2$ orbit 
under $f$.
We find $\log |F'(q_1)| = \log |f'(q_1) \cdot f'(q_2)| = \log 4|1-a|$.
So $\log |F'(q_1)| = \log 2$ only if $a = \frac12$ or $\frac32$, 
but neither parameter value corresponds to the required
renormalizable map.
\end{proof}

\section{ Liv\v{s}ic theorems for non-uniformly hyperbolic systems}

In this section we use
Young towers~\cite{Young1,Young2} to prove measurable Liv\v{s}ic
theorems for Lie group valued cocycles over non-uniformly hyperbolic systems.
In particular we are able to prove Liv\v{s}ic theorems for 
H\'enon maps  \cite{BC}. These maps take the form
$f(x,y)=(1-ax^2+y,\,bx),$ with $a\simeq 2$ and $b\simeq 0,$
where for a positive Lebesgue measure subset of parameter space $(a,b)$,
it is proved that these maps admit a
nontrvial attracting set with an ergodic Sinai-Ruelle-Bowen measure supported 
on it, see \cite{BY}. A Markov extension can be associated to such maps, as we describe below.
Our results are also applicable to other 
non-uniformly hyperbolic systems which can be shown to admit a Young tower, 
see for example \cite{ALV,Gou2,WY}.

We discuss diffeomorphisms for which there exists a stable foliation. 
For the non-uniformly expanding case drop all references to 
the stable foliation. We refer to Young's
original papers~\cite{Young1,Young2} and 
Baladi's book~\cite{Baladi} for more details.
In a non-uniformly hyperbolic system the unstable leaves
are not invariant under the return map to a reference set and this
introduces some complications to the analysis. 
The proof in this section is based on 
that of Theorem~\ref{young_tower} and Corollary~\ref{original_system}.

In this section we use notation that is commonly used for Young towers.
Let $T:X\to X$ be a $C^{1+\eps}$ diffeomorphism, where $X$ is a
compact manifold with metric $\rho_X$. 
Suppose there exists $\Lambda \subset X$ with a hyperbolic 
product structure ~\cite[Definition 1]{Young1}
 $\Lambda=\{(\cup \gamma^u)\cap(\cup \gamma^s)\,:\,\gamma^u\in\Gamma^u,\,\gamma^s\in\Gamma^s\}$, 
where $\Gamma^u,\Gamma^s$ are two families of $C^1$ disks in $X$ 
with the following properties: (i) 
disks  in $\Gamma^u$ are pairwise disjoint,
and the disks in $\Gamma^s$ are pairwise disjoint, 
(ii) every $\gamma^u\in\Gamma^u$ meets every $\gamma^s\in\Gamma^s$ in exactly
one point, 
(iii) there exists a lower bound on the  angle between $\gamma^u$ and $\gamma^s$ at the
point of intersection, and
(iv) each $\gamma^u\in\Gamma^u$ satisfies $m_{\gamma} (\gamma^u \cap \Lambda)>0,$
where 
$m_{\gamma}$ is the measure on $\gamma^u$ induced by the Riemannian structure of $X$.  

Under assumptions P1-P5 ~\cite[Section 1]{Young1}, Young constructs
a Markov extension (Young tower) $(F,\Delta)$ over $T:X\to X$ with
base $\Lambda$. The set $\Lambda$ is decomposed as $\Lambda=\cup_j \Lambda_j$
and there is a return function $R:\Lambda\to \N$, with constant value
$R_j$ on each $\Lambda_j$. Define $T^R (x)=T^{R(x)}(x)$.

$\Delta:=\{(x,l):x\in \Lambda;~l=0,1,\dots,R(x)-1\}$
where $\Delta_0=\{(x,0):x\in \Lambda\}=\Lambda$ in a natural
identification. 
Define $F:\Delta \to \Delta$ as 
\[
F(x,i) = \left\{ \begin{array}{ll}
(x,i+1) & \mbox{ if } x \in \Lambda_j \mbox{ and } i < R_j-1,\\
(T^{R_j}x,0) & \mbox{ if } x \in \Lambda_j \mbox{ and } i = R_j-1.
\end{array} \right.
\]
There is a projection $\pi:\Delta \to X$ given by $\pi(x,i) = T^i(x)$
and $\pi \circ F = T \circ \pi$.
The tower $(F,\Delta)$ is then reduced to an expanding map
$\bar{F}^R:\bar{\Delta}_0\To \bar{\Delta}_0$, where 
$\bar{\Delta}_0$ is the quotient of $\Delta_0$ under the equivalence
relation that two points are equivalent if and only if
they belong to the same local stable leaf $\gamma^s$.

We need the following properties.

(A1) There exists an $F$-invariant probability measure
$\nu$ with conditional measures $\{\nu^{\gamma}\}$ on
$\gamma^u\cap \Lambda$ leaves with  densities $\{\rho^{\gamma}\}$
with respect to $m_{\gamma}$
which  satisfy: $\frac{1}{C}\le \rho^{\gamma} \le C$ 
for some $C>0,$~\cite[Section 2]{Young2}.
The measure $\nu$ induces a $T$-invariant measure $\nu_X$ on $X$ defined 
by $\nu_X = \pi^* \nu$.

(A2) There is a countable partition $\P_0$ of $\Lambda=\Delta_0$ 
into elements $\{\Lambda_j\}$ together with a return function 
$R:\Lambda\to\mathbb{N}$, with $R\mid_{\Lambda_i}=R_i$. Moreover
$F^{R_i}$ maps $\Lambda_i$ bijectively onto 
$\Lambda$~\cite[Section 1.1]{Young2}. There is a countable partition $\P_\Delta$ with elements $(\Lambda_i,j)$ with $0 \leq j < R_i$ of $\Delta$.

(A3) There exists $K>0$ such that if $y\in \gamma^u(x):$
\[
\frac{1}{K}\le \frac{\|D^{u}(F^j x)\|}{\|D^{u}(F^j y)\|} \le K,
\]
for all $j=0,\dots,s(x,y)$~\cite[P4(b)]{Young2}. Here $D^{u}$ denotes the
derivative along unstable leaves $\gamma^{u}$, and $s(x,y)$ is the
first time $n$ for which $F^nx$ and $F^ny$ lie in different elements of $\P_0$.   

(A4) For $\gamma,\gamma'\in\Gamma^{u}$, if $\Theta:\gamma\cap\Lambda\To
\gamma' \cap \Lambda$ is defined by $\Theta(x)=\gamma^s(x)\cap \gamma',$
then $\Theta$ is absolutely continuous, and there exists a $C_1>0$
such that $$\frac{d(\Theta^{-1}_{*}\mu_{\gamma'})}{d\mu_{\gamma}}(x)\leq C_1,$$
for all $x\in\gamma^{u},$~\cite[P5(b)]{Young2}.

(A5) There exists $\lambda_u>1$ such that for each $x\in\Lambda$,  
$|D^{u}F^R(x)|\geq\lambda_u,$~\cite[Section 3.1]{Young2}.

(A6) There exists $\lambda_s<1$ such that for all $\gamma^s\in\Gamma^s$
and every $x,y$ in the same $\gamma^s$,
$\rho_X(F^jx,F^jy)\leq C\lambda^{j}_{s}$, cf. \cite[P3]{Young2}. 
Here $C>0$ is a uniform constant.

(A7)   Let ${\mathcal P}_n[x]$ denote the element of the partition
$ (T^{R})^{-n}{\mathcal P}_0$ that contains $x\in \Lambda $. For $y\in {\mathcal P}_n[x]$
 define $\tau(y)=R(y)+ R(T^R y)+\ldots + R((T^R)^{n-1 }y)$.
 Note that $\tau (y)=\tau (z)$ if $y,z\in {\mathcal P}_n [x]$.  Let $A_n$ be an element of
 the partition $ (T^{R})^{-n}{\mathcal P}_0$. Given $x_0\in \Lambda$ and 
 $x_{\tau_n}\in A_n$ with $T^{\tau_n}x_{\tau_n}=x_0$ 
define $x_i={T^{\tau_n -i}}x_{\tau_n}$. 
Let $\hat x=(x_0,x_{\tau_1}, x_{\tau_2},\ldots,x_{\tau_n},\ldots)$
be a point in  the natural extension of $T^{R}:\Lambda \to \Lambda$ with corresponding invariant measure 
$\hat \nu$. We assume either

\begin{enumerate}
\item  a one-dimensional   unstable manifold or
\item    for $\hat \nu$-a.e. $\hat x$ 
 there exists $C(\hat x)$ such that for all $z_{\tau_n},y_{\tau_n}\in
 A_n$, $0\le i\le \tau_n$, 
 \begin{eqnarray}\label{XXXX}
\rho_X(y_i,z_i)&\le& C(\hat x) \lambda_u^{i} \rho_X(y_0,z_0).
\end{eqnarray}
\end{enumerate}

\noindent {\bf Cocycle assumptions:}
Let $\phi\colon X\to G$ be H\"{o}lder of 
exponent $\alpha>0$, and define constants $\mu_u$ (as in 
Section~\ref{sec:groups}) and $\mu_s$ by:
\begin{eqnarray*}
\mu_u&:=&\lim\limits_{n\rightarrow\infty}\left(\sup\limits_{x\in
X}\|\Ad(\phi_n(x))\|\right)^{\frac{1}{n}},\\
\mu_s&:=&\lim\limits_{n\rightarrow\infty}\left(\sup\limits_{x\in
X}\|\Ad(\phi_n(x))^{-1}\|\right)^{-\frac{1}{n}}.
\end{eqnarray*}
where $\phi_n(x)=\phi (T^{n-1} x)\ldots \phi(x)$.
Assume  a partial hyperbolicity condition (\ref{PH})
on the group extension:
\begin{equation}\label{PH1}
 \lambda^{\alpha}_{s}<\mu_s\leq 1\leq  \mu_u < \tilde\lambda_{u}^{\alpha} \tag{PH}
\end{equation}
where $\tilde\lambda_u=\lambda_{u}^{\frac{1}{\mathcal{R}}}$.

We define:   
\[
\Theta_{PH}:=\max\biggl\{\frac{\log\mu_u}{\log\tilde\lambda_u},\,\frac{\log\mu_s}{\log\lambda_s}
\biggr\}.
\]
The assumption (\ref{PH1}) implies that $\Theta_{PH}<\alpha$.

\begin{theorem}\label{holder_tower}
Assume that $(T,X,\nu)$ is modelled by a tower over a base set 
$\Lambda\subset X$. Suppose in addition that $T$ has
a one-dimensional unstable direction or condition~(\ref{XXXX}).
Let $\phi: X\to G$ be H\"{o}lder of exponent $\alpha$ and suppose
condition (PH) holds. 
If  $\psi(T x)=\phi(x)\psi(x)$ $\nu_X$-a.e. for some 
measurable function $\psi: X\to G$, then 
$\psi\mid\Lambda$ is $\gamma$-H\"older for some $\gamma\in(0,1)$.
\end{theorem}

\begin{remark}
The same regularity is forced upon solutions $\psi$ to equation~\eqref{eq-cohomology1},
equation~\eqref{eq-cohomology2} or  equation~\eqref{eq-cohomology3}.
\end{remark}
\begin{remark}
Condition (PH) is automatic if $G$ is Abelian, compact or nilpotent.
\end{remark}

\begin{proof}

We start with a lemma tackling the stable direction.
\begin{lemma}
There exists $\psi'=\psi$ $\nu_X$-a.e.  and $\psi'$ is H\"{o}lder
when restricted to each $\gamma^s\in \Lambda$
(with uniform constant and exponent). 
\end{lemma}

\begin{proof}

Choose a version of $\psi$ and $\gamma^u\in \Lambda$ so that 
for $\nu_{\gamma^u}$-a.e. $\! z\in \gamma^u$, $\psi(Tz)=\varphi(z)\psi(z)$.  
For each $z\in \gamma^u$, each $x\in \gamma^s(z)\subset \Lambda$
define 
\[
\psi'(x)=\lim_{n\rightarrow \infty}\phi(x)^{-1}\ldots \phi(T^n x)^{-1}
\phi(T^n z)\ldots
\phi(z)\psi(z).
\]
By conditions (PH) and (A6), an argument similar to that of 
Section~\ref{sec:markov}
can be used to show that $\psi'$ restricted to each $\gamma^s(z)$
is uniformly H\"{o}lder. 
Furthermore  $\psi'(Tx)=\varphi(x)\psi'(x)$ for $\nu$-a.e. $x\in \Lambda$
and hence $\psi'=\alpha \psi$ $\nu$-a.e. for some constant group element
$\alpha$. As $\nu_X = \pi^*\nu$ the result follows.
\end{proof}
From now on we assume that $\psi$ has the properties specified in the lemma
above, namely $\psi$ restricted to each 
$\gamma^s\in \Lambda$ is uniformly H\"{o}lder.
Now we need only show that $\psi$ restricted to each $\gamma^u$ is
H\"{o}lder since the local product structure implies in this case 
that $\psi$ is H\"{o}lder on $\Lambda$.
In fact, to show that $\psi$ restricted to each $\gamma^u$ is H\"{o}lder
we need only show that there is a $\gamma^u\in \Lambda$ such
that $\psi$ restricted to $\gamma^u$ is H\"{o}lder, since the fact
that the holonomy is H\"{o}lder and $\psi$ restricted to each $\gamma^s$
is H\"{o}lder implies the result for all $\gamma^u\in \Lambda$.

Recall that $\P_0$ is the partition of $\Lambda$ into
$\{\Lambda_j \},$ and each $\Lambda_j$ contains whole stable leaves
$\gamma^s$. We will, to simplify notation, denote $(x,0)$ as $x$.
For $i \geq 1$ let
$\P_i = \bigvee_{j=0}^{i-1} (F^R)^{-j} \P_0$. Refine
the partition $\{ \P_i\}$ in the stable direction  by partitioning the stable manifolds
into leaves of length at most $2^{-i}$ to form a partition $\mathcal{Q}_i$
of $\Lambda$. Partition $\Lambda$ in such a way that if $A\subset B$, $A\in \mathcal{Q}_i,B\in \P_i$
then $\gamma^u\cap B=\gamma^u\cap A$ for each $\gamma^u$ such that 
$\gamma^u \cap A \neq \emptyset$.

The $\sigma$-algebra generated by $\bigvee_i\mathcal{Q}_i$ generates
the Borel $\sigma$-algebra on $\Lambda$. 
By the Martingale Density Theorem, given $\eta>0$ there
exists $n$ and an element of the partition  component  $A_{n} \in
\mathcal{Q}_n$ so that for some $x\in A_n$ 
\[
\frac{\nu\{ y \in A_n :d(\psi(x),\psi(y))<\eta\}}{\nu (A_n)}>1-\eta.
\]
For $x\in A_n$ define
 $\tau_n=R(x)+R(T^R x)+\ldots + R((T^R)^{n-1} x)$
and note $T^{\tau_n}(x)=(T^R)^n x$. Given $x_0\in\Lambda$,
 $x_{\tau_n}\in A_n$ with $T^{\tau_n}x_{\tau_n}=x_0$ define 
 $x_i$ by $T^{\tau_n-i} x_{\tau_n}=x_i$.
By (A1) and (A4) we may choose a portion of leaf
 $\gamma_n=\gamma^u\cap A_n$ and $x_{\tau_n}$  such that
 \[
\frac{m_{\gamma_n}\{ y \in \gamma_n\cap
  A_n:d(\psi(x_{\tau_n}),\psi(y))<\eta\}}{m_{\gamma_n} (\gamma_n)}
>1-\O(\eta).
\]

Then $\tilde\gamma_n:=(T^R)^n \gamma_n$ is  an unstable leaf
which crosses $\Lambda$ completely in the unstable direction.
 As a  consequence of (A3)
\begin{eqnarray}\label{tail}
\frac{m_{\tilde\gamma_n}\{ y_0 \in \tilde\gamma_n:
d(\psi(x_{\tau_n}),\psi(y_{\tau_n}))<\eta\}}
{m_{\tilde\gamma_n} (\tilde \gamma_n)}
>1-\O(\eta).
\end{eqnarray}

On  $\tilde\gamma_n$ define a  function 
$\Psi_n:\tilde\gamma_n  \To G$ by 
\[
\Psi_n (y_0)= \phi_{\tau_n} (y_{\tau_n}) [\phi_{\tau_n} (x_{\tau_n})]^{-1},
\]
where $\phi_{i}(x_j)=\phi(x_{j-i+1})\ldots\phi(x_{i+j})$.

Take points $z_0,w_0\in\tilde\gamma_n$. Then by the cohomological equation 
\begin{eqnarray*}
\psi(z_0)&=&\phi_{\tau_n} (z_{\tau_n}) \psi(z_{\tau_n})\\
&=&\Psi_n (z_0)[\phi_{\tau_n} (x_{\tau_n})]\psi(x_{\tau_n})\psi(x_{\tau_n})^{-1}\psi (z_{\tau_n})\\
&=& \Psi_n (z_0)\psi(x_0)\psi(x_{\tau_n})^{-1}\psi (z_{\tau_n}).
\end{eqnarray*}

By the  the right-invariance of the metric and the triangle inequality we have
\begin{eqnarray*}
d(\psi(z_0),\psi(w_0))&\leqslant&d(\Psi_n(z_0)\psi(x_0)\psi(x_{\tau_n})^{-1}\psi (z_{\tau_n}),
\Psi_n(z_0)\psi(x_0))\\
&+& d(\Psi_n(z_0)\psi(x_0)\psi(x_{\tau_n})^{-1}\psi (w_{\tau_n}),
\Psi_n(z_0)\psi(x_0))   \\
&+&d(\Psi_n(z_0)\psi(x_0)\psi(x_{\tau_n})^{-1}\psi (w_{\tau_n}),
\Psi_n(w_0)\psi(x_0)\psi(x_{\tau_n})^{-1}\psi (w_{\tau_n})).
\end{eqnarray*}

We claim that  $\Psi_n$ is H\"{o}lder on $\tilde \gamma_n$ with 
uniform H\"{o}lder constant and exponent (the uniformity is 
over $n$ in the construction).

We calculate
\begin{eqnarray*}
\lefteqn{
d(\phi_{\tau_n}(z_{\tau_n})\phi_{\tau_n}(x_{\tau_n})^{-1},
\phi_{\tau_n}(w_{\tau_n})\phi_{\tau_n}(x_{\tau_n})^{-1})} \qquad \\ 
&=& d(\phi_{\tau_n}(z_{\tau_n}),\phi_{\tau_n}(w_{\tau_n})) \\
&\leqslant&\slims_{i=0}^{\tau_n-1}\|\Ad(\phi_i(z_i))\|d(\phi(z_{i+1}),\phi(w_{i+1}))\\
&\leqslant&\slims_{i=0}^{\tau_n-1}C (\mu_u  )^{i+1} \lambda_u^{-(i+1)\alpha}\rho_X(z_0,w_0)^{\alpha},
\end{eqnarray*}
where in passing from the second to third line, we use condition~(\ref{XXXX}).
Equivalently for the unstable direction we could have used 
the existence of the positive Lyapunov exponent $\tilde\lambda_u$, 
together with bounded distortion as in the proof of Theorem~\ref{main}. 
The series converges uniformly because of condition (PH).
Equation~(\ref{tail}) gives 

\begin{eqnarray}
\frac{m_{\tilde\gamma_n}\{ y_0 \in \tilde\gamma_n:
d(\psi(x_{0}),\psi(y_{0}))<C\rho_X(x_0,y_0)^{\alpha}\}}{m_{\tilde\gamma_n} (\tilde \gamma_n)}
>1-\O(\eta).
\end{eqnarray}  
Recall that $\psi$ restricted to each stable leaf in $\Lambda$ is uniformly H\"{o}lder.
The holonomy map along stable leaves is absolutely continuous (A4) and the density
of $\nu$ with respect to Lebesgue is bounded away from zero and above by (A1). 
Hence \cite[Proposition 19.1.1]{KH}
implies that 
 
\begin{eqnarray}
\frac{\nu\{ (x,y) \in \Lambda\times \Lambda:
d(\psi(x),\psi(y))<\tilde C \rho_X(x,y)^{\alpha}\}}{\nu\times\nu(\Lambda\times\Lambda)}
>1-\O(\eta).
\end{eqnarray} 

Since $\nu_X = \pi^*\nu$ and $\eta$ is arbitrary, it follows that
$\psi|_{\Lambda}$ has a H\"{o}lder version, thus proving Theorem~\ref{holder_tower}.
\end{proof}

\section{Appendix}
Suppose $\psi: X\to G $ is a measurable function from a metric
measure space into a connected finite-dimensional matrix  Lie group $G$ endowed with a right
invariant metric $d_G$. Let 
$\pi_{i,j}:G \to \R$, $1\le i,j\le 2d$  be the local coordinate
chart  functions.
 The proof of the proposition  below clearly generalizes to  any finite number of
 real-valued measurable functions
and hence establishes the Martingale Convergence Theorem, since
continuity is a local property.

\begin{proposition}
Suppose $(X,\mu)$ is a probability space and  $\{{\mathcal P}_n\}$ is an increasing sequence
of partitions of $X$ and let ${\mathcal P}_n [x]$
denote the partition element of ${\mathcal P}_n$ which
contains $x\in X$. Suppose the Borel $\sigma$-algebra is generated
by $\bigvee_n {\mathcal P}_n$. Let  $\phi: X\to \R$  be $\mu$-measurable and $\eta>0$. For $\mu$-a.e. $x\in X$, 
\begin{eqnarray}\label{good_points}
\lim_{n\to\infty}\frac{\mu \{y\in {\mathcal P}_n [x]: 
d(\phi(x),\phi(y))<\eta\}}{\mu({\mathcal P}_n)}>1-\eta.
\end{eqnarray}
\end{proposition}
\begin{proof}
First suppose $\phi\in L^1 (\mu)$.
Let ${\mathcal F}_n$ denote the $\sigma$-algebra generated by the partition ${\mathcal P}_n$.
Then 
\[
\lim_{n\to \infty} \E[\phi|{\mathcal F}_n](x)=\phi (x),\qquad
\mu-a.e. 
\]
by~\cite[Corollary 5.22]{Breiman}. Note  $\E[\phi|{\mathcal
  F}_n](x)$ is constant on ${\mathcal P}_n$.
Choose a sequence    $\{\delta_i\}$ such that $\sum_i {\delta_i}<\infty$. 
Given $\delta_i>0$  take $N_i$ sufficiently large that 
 $d(\E[\phi|{\mathcal F}_n](x),\phi (x))< \eta$ except for a 
 set of measure at most ${\delta_i}\eta^2$ for all $n\ge N_i$. For all
 $N_i$, the  union $U_{i}$ 
 of the set of atoms $A\in  {\mathcal P}_{N_i}$ for which 
 \[
\frac{\mu \{y\in A: 
d(\E[\phi|{\mathcal F}_{N_i}](y),\phi(y))>\eta\}}{\mu(A)} > \eta,
\]
 satisfies $\mu(U_i)<\eta^2\delta_i$. By the Borel-Cantelli Lemma, $\mu$-a.e. 
$x\in X$ lies in only finitely many $U_i$. 
Finally to remove the assumption that $\phi$ is integrable note that given
 $\eps>0$  there exists an integrable function $\psi$ such 
that $\psi(x)=\phi(x)$
 except for a set  of measure at most $\eps$. 
An argument using approximating functions and the Borel-Cantelli Lemma 
gives the  same result for measurable $\phi$. 
 \end{proof}

\medskip
\noindent
Henk Bruin\\
Mathematics and Statistics\\
University of Surrey\\
Guildford, Surrey, GU2 7XH\\
UK\\
\texttt{h.bruin@surrey.ac.uk}\\
\texttt{http://www.maths.surrey.ac.uk/showstaff?H.Bruin}

\medskip
\noindent
Mark Holland\\
Mathematics and Statistics\\
University of Surrey\\
Guildford, Surrey, GU2 7XH\\
UK\\
\texttt{mark.holland@surrey.ac.uk}\\
\texttt{http://www.maths.surrey.ac.uk/showstaff?M.Holland}

\medskip
\noindent
Matt Nicol\\
Mathematics \\
University of Houston \\
Houston TX 77204-3008\\
USA\\
\texttt{nicol@math.uh.edu}\\
\texttt{http://www.math.uh.edu/}

\end{document}